\documentclass[leqno,12pt]{amsart} 
\usepackage{amssymb} 
\usepackage{amsmath}
\usepackage{amsthm}
\usepackage{bm}
\usepackage{graphicx,psfrag,wrapfig}
\theoremstyle{plain} %
\newtheorem{theorem}{\indent\sc Theorem}[section] %
\newtheorem{lemma}[theorem]{\indent\sc Lemma}

\newtheorem{proposition}[theorem]{\indent\sc Proposition}

\theoremstyle{definition} %
\newtheorem{definition}[theorem]{\indent\sc Definition}
\newtheorem{remark}[theorem]{\indent\sc Remark}
\newtheorem{example}[theorem]{\indent\sc Example}

\newtheorem{fact}[theorem]{\indent\sc Fact}
\newcommand{\affine}{\mathbb{C}}
\newcommand{\Hproduct}[2]{\langle #1 , #2 \rangle}
\newcommand{\tangentvector}[1]{\partial / \partial #1}
\newcommand{\norm}[1]{\left|\!\left|#1\right|\!\right|_{\infty}}
\newcommand{\moduli}{\mathcal{M}}
\newcommand{\elliptic}{\mathbb{C} / \Lambda}
\newcommand{\degree}[1]{\mathrm{deg}\, (#1)}
\newcommand{\Dnorm}[3]{\left|\!\left| #1 \right|\!\right|_{\mathcal{C}^{#2}(#3)}}

\begin{document}

\title[A Packing Problem for Holomorphic Curves]{A Packing Problem for Holomorphic Curves} 

\author[Masaki Tsukamoto]{Masaki Tsukamoto}




\maketitle

\begin{abstract}
We propose a new approach to the value distribution theory of entire holomorphic curves.
We define a \textit{packing density} of an entire holomorphic curve, and show that 
it has various non-trivial properties.
We prove a \textit{gap theorem} for holomorphic maps from elliptic curves to the complex projective space,
and study the relation between theta functions and our packing problem.
Applying the Nevanlinna theory,
 we investigate the packing densities of entire holomorphic curves 
in the complement of hyperplanes.

\end{abstract}

\section{Main results}
\subsection{Introduction}
Since R. Nevanlinna discovered his celebrated theory on meromorphic functions [N], thousands of researchers have studied the value distribution theory of meromorphic functions and, more generally,
entire holomorphic curves in complex manifolds.
 
This paper is a new approach to the value distribution theory.
We define and study a \textit{packing problem} for entire holomorphic curves.

 \textit{Packing} is usually a notion in discrete geometry. 
For example, the Kepler conjecture on the sphere packing in the three dimensional Euclidean space is 
very famous.
In this paper we define a \textit{packing density} of an entire holomorphic curve and investigate its behavior. (To be precise, we consider only entire holomorphic curves with bounded derivative.)
 In particular, we study entire holomorphic curves in the complex projective space and prove that 
 their packing densities have a \textit{non-trivial} upper bound.

It is well-known that lattices in the Euclidean space play an important role in the sphere packing problem [CS]. They give the systematic way to pack spheres, and such lattice packings sometimes become very efficient configuration.
 In a similar way, holomorphic maps from elliptic curves play an important role in the packing problem 
 for entire holomorphic curves.
They sometimes produce good lower bounds for the packing problem.
 We apply these lower bounds to the projective embeddings of elliptic curves by means of 
theta functions, and show that they \textit{asymptotically} give the best packing.
In addition, as a different application of the lower bounds, we can prove a certain \textit{gap theorem} for holomorphic maps from elliptic curves to the complex projective space.

In the final section of the paper, we study the relation between the packing problem and the usual value distribution theory.
As is well-known,
 if a meromorphic function in the complex plane omits more than two values, it reduces to a constant function (Picard's theorem). 
In this paper, we investigate meromorphic functions which  omit exactly two values, 
and show that their behavior in the packing problem is  quite different from 
that of general meromorphic functions.
We also prove a more general result on the packing densities of entire holomorphic curves in the complement of hyperplanes in the complex projective space.

\subsection{Packing density and holomorphic capacity}
To begin with, we fix conventions on Hermitian manifolds. The basic reference is [KN].
Let $X$ be a complex manifold and $g$ be a Hermitian metric on $X$, i.e., $g$ is a Riemannian metric
on X compatible with the complex structure $J$:
\[g(Ju,Jv) = g(u,v) \quad \text{ for any vector fields $u$ and $v$.}\]
We extend $g$ to a complex symmetric bilinear form on $T^\affine X := TX \otimes_\mathbb{R} \affine$,
and also denote it by $g$.
From $g$, we can define the Hermitian inner product $\Hproduct{\cdot}{\cdot}$ and the norm $|\cdot|$ on $T^\affine X$ by
\[\Hproduct{u}{v} := g(u, \bar{v}), \quad |u| := \sqrt{\Hproduct{u}{u}},
 \quad \text{for all } u, v \in \Gamma(T^\affine X). \]
Here $\bar{v}$ is the complex conjugate of $v$.

The fundamental 2-form $\omega$ is defined by 
\[\omega(u, v) := g(Ju, v) \quad \text{for all } u, v \in \Gamma(T^\affine X). \]
(This is different from the convention in [KN].)

In terms of a local coordinate system $(z_1, z_2, \dots, z_n)$ on $X$, ($n := \mathrm{dim}_\affine X$),
the fundamental 2-form $\omega$ is expressed by
\[\omega = \sqrt{-1} \, \sum g_{i\bar{j}} dz_i \wedge d\bar{z}_j \quad (g_{i\bar{j}} := 
\Hproduct{\tangentvector{z_i}}{\tangentvector{z_j}}) .\]

Let $f: \affine \to X$ be a holomorphic curve, i.e., a holomorphic map from the complex plane $\affine$ to $X$.
Let $z = x + y \sqrt{-1}$ be the natural coordinate on the complex plane.
We define a norm of the differential $df : T^{\affine} \affine \to T^{\affine} X$ by setting 
\begin{equation}
|df|(z) := \sqrt{2} \, |df(\tangentvector{z})|.  \label{def:norm of the differential}
\end{equation}
Here $\tangentvector{z} = \frac{1}{2} \,(\tangentvector{x} - \sqrt{-1} \tangentvector{y})$, and 
the factor $\sqrt{2}$ in (\ref{def:norm of the differential}) 
comes from $|\tangentvector{z}| = 1/\sqrt{2}$. (The complex plane is equipped with the usual Euclidean metric.) The norm $|df|$ satisfies 
\begin{equation}
f^* \omega = |df|^2 \, \frac{\sqrt{-1}}{2} \, dz \wedge d\bar{z} = |df|^2 dx \wedge dy.
\end{equation}

Next we define the moduli space of holomorphic curves by
\begin{equation}
\moduli (X, \omega) := \{ f: \affine \to X | \text{ $f$ is holomorphic and } |df|(z) \leq 1
\text{ for all z $\in \affine$} \}. \label{def:moduli}
\end{equation}
For a holomorphic curve $f \in \moduli (X, \omega)$, we define the \textit{packing density} 
$\rho(f)$ by setting 
\begin{equation}
\rho(f) :=  \limsup_{R \to \infty} \frac{1}{\pi R^2}\int_{|z|\leq R} f^*\omega
      \, = \, \limsup_{R \to \infty} \frac{1}{\pi R^2}\int_{|z|\leq R} |df|^2 \, dxdy. 
\label{def:density}
\end{equation}
This satisfies 
\begin{equation*}
0 \leq \rho(f)  \leq 1.
\end{equation*}
The integration of $f^* \omega$ is the energy functional. 
Hence if $\rho(f)$ is close to $1$, the energy of $f$ is \textit{densely packed} in the complex plane.
In other words, $\rho(f)$ evaluates the efficiency of the energy distribution 
of a holomorphic curve $f$. 
This is the reason why we call $\rho (f)$ ``packing density''.

Using the packing density $\rho(f)$, we define the \textit{holomorphic capacity} $\rho(X, \omega)$
by setting 
\begin{equation}
\rho(X, \omega) := \sup_{f \in \moduli (X, \omega)} \rho(f).  \label{def:capacity}
\end{equation}
This satisfies the following.
\begin{equation*}
0 \leq \rho(X, \omega) \leq 1.
\end{equation*}

Then we can define a packing problem.
The packing problem for holomorphic curves is the problem of determining, or estimating, the value of 
the holomorphic capacity $\rho (X, \omega)$.
We will often abbreviate $\moduli (X,\omega)$ and $\rho (X, \omega)$ to $\moduli(X)$ and $\rho(X)$
when it causes no confusion.
\begin{example}
Let $(z_1, z_2, \cdots, z_n)$ be the natural coordinate system on $\affine^n$. 
The Euclidean metric and its fundamental 2-form on $\affine^n$ are expressed by
\begin{equation*}
ds^2 = \sum_{i = 1}^n dz_i d\bar{z}_i,  
\quad \omega = \frac{\sqrt{-1}}{2}\sum_{i=1}^n dz_i \wedge d\bar{z}_i.
\end{equation*}
Let $f:\affine \to \affine^n$ be the natural inclusion:
$f(z) := (z, 0, 0, \cdots, 0)$.
It is obvious that $|df| \equiv 1$. Hence $f \in \moduli(\affine^n, \omega)$ and $\rho(f) = 1$. 
Therefore 
\begin{equation*}
\rho(\affine^n, \omega) = 1.
\end{equation*}

In the same manner, if $X$ is a complex torus with the Euclidean metric induced by the universal covering, we have
\begin{equation*}
\rho(X) = 1.
\end{equation*}
\end{example}
\begin{example}
Let $\Delta = \{ z \in \affine \,| \, |z| < 1 \}$ be the unit disk and $g$ be an arbitrary Hermitian metric on $\Delta$.
All holomorphic maps from $\affine$ to $\Delta$ are constant maps by Liouville's theorem.
Therefore we have
\begin{equation*}
\rho (\Delta, \omega) = 0. 
\end{equation*}

In the same manner, if $X$ is a compact Riemann surface with a genus $\geq 2$, we have
\begin{equation*}
\rho (X, \omega) = 0 \quad \text{for any Hermitian metric on $X$.}
\end{equation*}
\end{example}
The above two examples are trivial extremal cases. Our main concern is the case
of the complex projective space $\affine P^n$ with the Fubini-Study metric, ($n \geq 1$).
To state the results clearly, we explicitly define the Fubini-Study metric as follows.

Let $[z_0:z_1:\cdots:z_n]$ be the homogeneous coordinate in $\affine P^n$ and set
\[U_0 := \{ \, [1: z_1 : z_2 : \cdots : z_n ] \in \affine P^n |\,
 (z_1, z_2, \cdots , z_n) \in \affine ^n \} . \]
We define the Fubini-Study metric form $\omega_{FS}$ on $U_0$ by 
\begin{equation}
  \begin{split}
\omega_{FS} &:= \frac{\sqrt{-1}}{2 \pi} \partial \bar{\partial}
                \log \left( 1 + \sum_{i = 1}^n |z_i|^2 \right),\\
            &= \frac{\sqrt{-1}}{2 \pi} \sum_{1\leq i,j \leq n} \frac{\partial^2}{\partial z_i \partial \bar{z}_j}
               \log \left( 1 + \sum_{k=1}^n |z_k|^2 \right) \, dz_i \wedge d\bar{z}_j
\quad \text{on $U_0$}.
  \end{split}    \label{def:Fubini-Study}
\end{equation}
It is standard that this 2-form $\omega_{FS}$ smoothly extends over $\affine P^n$ and define the Fubini-Study metric. This is normalized so that
\begin{equation}
\int_{\affine P^1} \omega_{FS} = 1
 \quad \text{for $ \, \affine P^1 := \{ \, [z_0: z_1: 0: \cdots : 0] \in \affine P^n \}$ }.
\label{normalization}
\end{equation}
The first result is the following.
\begin{theorem}\label{thm:non-trivial}
\[  0 \, \lvertneqq  \, \rho (\affine P^n, \omega_{FS}) \,  \lvertneqq 1.  \] 
\end{theorem}
This result means that the holomorphic capacity is a non-trivial object.
We will usually abbreviate $\rho (\affine P^n, \omega_{FS})$ to $\rho (\affine P^n)$.

Theorem \ref{thm:non-trivial} does \textit{not} give an \textit{effective} upper bound for $\rho(\affine P^n)$. We investigate the explicit estimate for $\rho (\affine P^1)$ in the next theorem. 
\begin{theorem} \label{thm:effective}
\begin{equation*}
 \rho(\affine P^1)  \leq  1 - 10^{-100}.  
\end{equation*}
\end{theorem}
The above value, $1- 10^{-100}$, itself has no importance.
The important point is that it is an \textit{explicit} number.

Next we study behavior of $\rho(\affine P^n)$ as $n$ goes to infinity.
The natural inclusion $\affine P^n = \{ \, [z_0: z_1: \cdots : z_n: 0] \in \affine P^{n+1} \} \hookrightarrow \affine P^{n+1}$ is a holomorphic isometric imbedding. Hence we can consider
\begin{equation*}
\moduli (\affine P^1) \subset \moduli (\affine P^{2}) \subset \moduli (\affine P^{3})
\subset \cdots  \subset \moduli (\affine P^{n}) \subset \moduli (\affine P^{n+1}) \subset \cdots.
\end{equation*}
It results that 
\begin{equation}
0 \lvertneqq \rho(\affine P^1) \leq \rho(\affine P^2) \leq \rho(\affine P^3) \leq \cdots
\leq \rho(\affine P^n) \leq \rho(\affine P^{n+1}) \leq \cdots \lvertneqq 1.  \label{ascend}
\end{equation}
The following theorem determines the limit of this sequence.
\begin{theorem}\label{thm:limit}
\[ \lim_{n \to \infty} \rho(\affine P^n) \, = \, 1. \] 
\end{theorem}
We will prove this theorem by using theta functions.

\subsection{Application to holomorphic maps from elliptic curves}
As an application of holomorphic capacity, we will prove a ``\textit{gap theorem}''
for elliptic functions and, more generally, holomorphic maps from elliptic curves to $\affine P^n$.

Let $\elliptic$ be a elliptic curve. Here $\Lambda$ is a lattice in $\affine$ 
with $\mathrm{rank} \, \Lambda = 2 $. We give $\elliptic$ the Euclidean metric induced by 
the universal covering $\affine$.
Let $f: \elliptic \to \affine P^n$ be a holomorphic map. 
(Here we don't consider any restriction on the norm $|df|$.)
 We define the degree of $f$, $\degree{f}$, by setting
\begin{equation}
\degree{f} := \int_{\elliptic} f^* \omega_{FS} = \int_{\elliptic} |df|^2 dxdy.
   \label{def:deg}
\end{equation}
From the normalization of the Fubini-Study metric in (\ref{normalization}), $\degree{f}$ is a 
non-negative integer and a homological invariant.
From (\ref{def:deg}), we get the trivial estimate:
\begin{equation}
\norm{df}^2 \geq \frac{\degree{f}}{\mathrm{vol} \, (\elliptic)}. \label{trivial}  
\end{equation}
Here $\norm{df} := \sup_{z \in \elliptic} |df|(z)$, and $\mathrm{vol} \, (\elliptic)$ is the volume 
of $\elliptic$ defined by the Euclidean metric.

The following result shows that this is \textit{not} a best estimate.
\begin{theorem}\label{thm:gap}
For any holomorphic map $f: \elliptic \to \affine P^n$, we have 
\[ \norm{df}^2 \geq \frac{1}{\rho(\affine P^n)} \frac{\degree{f}}{\mathrm{vol} \, (\elliptic)}. \]
\end{theorem}
From Theorem \ref{thm:non-trivial}, 
\[ \frac{1}{\rho(\affine P^n)} \gvertneqq 1. \]
Hence there exists a certain \textit{gap} between the trivial estimate (\ref{trivial})
and Theorem \ref{thm:gap}. The point is that $1/ \rho(\affine P^n)$ is the universal constant which does not
depend on any lattice $\Lambda$ nor any holomorphic map $f$.


\subsection{Holomorphic capacity of the complement of hyperplanes}

For a meromorphic function $f(z)$ in the complex plane, a point $w \in \affine P^1$ is called 
a lacunary value if $f(z) \neq w$ for all $z$.
The famous Picard's theorem states that a non-constant meromorphic function in the complex plane has at most two lacunary values.
Hence ``2'' is the critical value of the number of lacunary values. 
We study the holomorphic capacity of this borderline case:
\begin{theorem}\label{thm:Nevanlinna}
Let $w_1$ and $w_2$ be two distinct points in $\affine P^1$. Then 
\[ \rho (\affine P^1 \setminus \{ w_1, w_2 \} ,\, \omega_{FS}) = 0 . \]
Here we use the Fubini-Study metric restricted to $\affine P^1 \setminus \{ w_1, w_2 \}$
as the metric on it.
\end{theorem}

This result can be generalized to the higher dimensional complex projective space as follows.
First we fix the notion of ``linearly independence'' of hyperplanes in the complex projective space.

Let $H_0, H_1, \cdots , H_{n}$ be the $n+1$ hyperplanes in $\affine P^n$ defined by
\begin{equation}\label{defining equations of hyperplanes}
 H_i:\quad  \sum_{j = 0}^n a_{i j} z_j  = 0, \quad (0 \leq i \leq n).  
\end{equation}
Here $[ z_0: z_1: \cdots : z_n ]$ is the homogeneous coordinate in $\affine P^n$.
Let $A := (a_{ij})_{0 \leq i, j \leq n}$ be the coefficients matrix. 
 $H_0, H_1, \cdots H_n$ are said to be linearly independent if $A$ is a regular matrix.
\begin{theorem}\label{thm:Nevanlinna 2}
Let $H_0, H_1, \cdots, H_n$ be $n+1$ linearly independent hyperplanes in $\affine P^n$.
Then 
\[ \rho (\affine P^n \setminus (H_0 \cup H_1 \cup \cdots \cup H_n), \, \omega_{FS}) = 0 .\]
Here we use the Fubini-Study metric restricted to 
$\affine P^n\setminus (H_0 \cup H_1 \cup \cdots \cup H_n)$ as the metric on it.
\end{theorem}

This theorem makes a sharp contrast with Theorem \ref{thm:non-trivial} and Theorem \ref{thm:limit}.

\begin{remark}
The moduli space $\moduli(\affine P^n \setminus (H_0 \cup H_1 \cup \cdots \cup H_n),\, \omega_{FS})$ has
non-constant holomorphic curves.
For example, $\exp : \affine \to \affine P^1 \setminus \{0, \infty \}$ satisfies
\[ |d \exp |(z) = \frac{1}{\sqrt{\pi}}\frac{e^x}{1+e^{2x}} < 1, \quad (x = \mathrm{Re}\,z). \]
Hence it is an element of $\moduli(\affine P^1 \setminus \{0, \infty\},\, \omega_{FS})$.
\end{remark}

\subsection{Organization of the paper}
In Section 2, we develop general theory of holomorphic capacity and prove Theorem \ref{thm:gap}.
 We show Theorem \ref{thm:non-trivial}
in Section 3 by applying the results of Section 2. 
We prove Theorem \ref{thm:limit} in Section 4 by using theta functions and the result of Section 2.
We prove Theorem \ref{thm:effective} in Section 5. 
Theorem \ref{thm:Nevanlinna} is a special case of Theorem \ref{thm:Nevanlinna 2},
 and we prove it in Section 6 by using the Nevanlinna theory.

Section 5 is independent of Section 3 and 4. Section 6 is logically independent of 
all other sections. (But its meaning in the packing problem is underpinned by other results.)

\section{General Theory}
In this section we study general properties of holomorphic capacity. 
We discuss its scaling invariance in Section 2.1.
 We study upper bounds for holomorphic capacity in Section 2.2,
 and we study lower bounds in Section 2.3.
In this section $X$ is a complex manifold with a Hermitian metric $g$ and its fundamental 
2-form $\omega$. 

\subsection{Scaling invariance}
Let $f:\affine \to X$ be a holomorphic map. 
Suppose that there exists a positive number $m \lvertneqq \infty$ such that 
\begin{equation}
|df|(z)  \leq m  \quad \text{for all $z \in \affine$}.
\end{equation}
We define the holomorphic map $\Hat{f}: \affine \to X$ by setting $\Hat{f}(z) := f(z/m)$.
Then 
\begin{equation}
|d\Hat{f}|(z) = \frac{1}{m} |df|(z/m) \leq 1. \label{scaling}
\end{equation}
Hence $\Hat{f} \in \moduli (X, \omega)$.
\begin{lemma} \label{lem:scaling}
\[ \limsup_{R \to \infty} \frac{1}{m^2 \pi R^2} \int_{|z| \leq R} |df|^2 dxdy = \rho(\Hat{f}). \]
\end{lemma}
\begin{proof}
From (\ref{scaling}) 
\begin{equation*}
\frac{1}{m^2 \pi R^2} \int_{|z| \leq R} |df|^2(z) \, \frac{\sqrt{-1}}{2} dz \wedge d\bar{z}
= \frac{1}{\pi R^2} \int_{|z| \leq R} |d\Hat{f}|^2(mz) \, \frac{\sqrt{-1}}{2}dz \wedge d\bar{z}.
\end{equation*}
Introducing the new variables $w := mz$ and $r := mR$, we get
\begin{equation*}
\frac{1}{m^2 \pi R^2} \int_{|z| \leq R} |df|^2(z) \, \frac{\sqrt{-1}}{2} dz \wedge d\bar{z}
= \frac{1}{\pi r^2}\int_{|w| \leq r} |d\Hat{f}|^2(w) \, \frac{\sqrt{-1}}{2} dw \wedge d\bar{w}.
\end{equation*}
Hence we get the conclusion.
\begin{equation*}
   \begin{split}
\limsup_{R\to \infty}
\frac{1}{m^2 \pi R^2} \int_{|z| \leq R} |df|^2(z) \, \frac{\sqrt{-1}}{2} dz \wedge d\bar{z}  &=
\limsup_{r \to \infty} 
\frac{1}{\pi r^2} \int_{|w| \leq r} |d\Hat{f}|^2(w) \, \frac{\sqrt{-1}}{2} dw \wedge d\bar{w}, \\
&= \rho(\Hat{f}).
   \end{split}
\end{equation*}
\end{proof}
The following result states that holomorphic capacity is invariant under a scale change of 
a Hermitian metric.
\begin{proposition} \label{prop:scaling}
For any positive number $c$, we have
\[ \rho(X,\, c \, \omega)  =  \rho(X, \, \omega).  \]
$\rho(X, \, c \, \omega)$ is
the holomorphic capacity defined by using $c \, g$ as the metric on $X$.
\end{proposition}
\begin{proof}
Set 
\[|u|' := \sqrt{c} \, |u|  \quad \text{for all $u \in \Gamma (T^{\affine} X)$}. \]
This is the norm induced by $c  g$. Let $f$ be an element of $\moduli (X, \, c\, \omega)$. 
Then $|df|' \leq 1$.

Set $\Hat{f}(z) := f(\sqrt{c}\, z)$, then 
\[ |d\Hat{f}|(z) = \sqrt{c} \, |df|(\sqrt{c} \, z) = |df|'(\sqrt{c} \, z) \leq 1. \]
This means that $\Hat{f}$ is an element of $\moduli (X, \, \omega)$.
Hence there exists the following one-to-one correspondence:
\[ \moduli (X, c\, \omega) \ni f(z)  \quad \longleftrightarrow \quad
 \Hat{f}(z) = f(\sqrt{c} \, z) \in \moduli (X, \omega). \]
From Lemma \ref{lem:scaling} with $m = 1/\sqrt{c}$,
\[ \limsup_{R \to \infty} \frac{1}{\pi R^2} \int_{|z| \leq R} |df|'^2 dxdy  =
   \limsup_{R \to \infty} \frac{1}{\pi R^2} \int_{|z| \leq R} |d\Hat{f}|^2 dxdy. \]
It follows that $\rho(X,\, c \, \omega)  =  \rho(X, \, \omega)$.
\end{proof}

\subsection{Upper bounds for holomorphic capacity} \label{subsection:upper bound}
In this section we study upper bounds for holomorphic capacity. 
First we will introduce a key notion. For a positive number $r$, let $\Delta(r)$ be the open disk 
of radius $r$ centered at the origin in the complex plane:
\[ \Delta (r) := \{ z \in \affine| \, \lvert z \lvert <  r \}. \]
We use the natural Euclidean metric as the metric on $\Delta(r)$.
\begin{definition}
\footnote{The word ``WFL'' is inspired by the word ``WFF'' in [DK, Definition (3.2.2)].}
A Hermitian manifold $X$ is $r$-WFL (without flat lines) if there is no
holomorphic isometric immersion from $\Delta(r)$ to $X$. 
$X$ is WFL if $X$ is $r$-WFL for all positive numbers $r$.

Here, a holomorphic isometric immersion from $\Delta(r)$ to $X$ is a holomorphic map
$f: \Delta(r) \to X$ with $|df|(z) = 1$ for all $z \in \Delta(r)$.
\end{definition}
\begin{example}
The complex projective line $\affine P^1$ is WFL.
\begin{proof}
Suppose that there exists a holomorphic isometric immersion from $\Delta(r)$ to $\affine P^1$ for some positive number $r$.
Since $\mathrm{dim}_{\affine} \Delta(r) = \mathrm{dim}_{\affine} \affine P^1$, this means that $\Delta(r)$ is locally isometric to $\affine P^1$. Because $\Delta(r)$ is flat and $\affine P^1$ has a positive constant curvature, it is impossible.
\end{proof}
In Section 3 we prove that all $\affine P^n$ are also WFL.
\end{example}

The following is the main result of this section.
\begin{theorem} \label{thm:upper-bound}
Let $X$ be a compact Hermitian manifold and suppose 
$X$ is $r$-WFL for some positive number $r$. Then 
\[ \rho(X) \lvertneqq 1.  \]
\end{theorem}
Before proving the theorem, we need a preliminary result.
To begin with, note that we can suppose $X$ is $1/2$-WFL without loss of generality by using a scale change and Proposition \ref{prop:scaling} if $X$ is $r$-WFL for a positive number $r$. (This is just for simplicity.)

The key proposition is the following.
\begin{proposition} \label{prop:key}
Let $X$ be a $1/2$-WFL compact Hermitian manifold and $K$ be a unit square in the complex plane $\affine$. 
Then, there exists a constant $c(K) \lvertneqq 1$ such that 
\[ \int_{K} |df|^2 \, dxdy  \,  \leq c(K)  \quad  \text{for all $f \in \moduli (X)$}. \]
\end{proposition}
\begin{proof}
Since $\mathrm{vol} (K) =1$ and $|df| \leq 1$, it is trivial that $\int_{K} |df|^2 \, dxdy  \,  \leq 1$.
Hence the proposition states that this trivial estimate can be improved.
 
First, note that a unit square contains a open disk of radius $1/2$.

Suppose the proposition is false. 
Then we have a sequence $\{ f_n \}_{n\geq 1}$ in $\moduli (X)$ such that 
\begin{equation}
\int_{K} |df_n|^2 \, dxdy  \to 1 , \quad  (n \to \infty).   \label{limit}
\end{equation}
Because $|df_n| \leq 1$ and $X$ is compact, we can apply 
Arzela-Ascoli's theorem, and get a continuous map $f: \affine \to X$ such that an appropriate 
subsequence of $\{ f_n \}_{n\geq 1}$ converges to $f$ in the sense of compact uniform convergence.
Since each $f_n$ is holomorphic, $f$ is a holomorphic map and, if we take a subsequence, 
\[ |df_n| \to |df|, \quad (n \to \infty) ,\]
in the compact uniform topology. 
From $|df_n| \leq 1$, we get $|df| \leq 1$. 
From the assumption (\ref{limit}),
\[ \int_{K} |df|^2 \, dxdy  =  1. \]
Since $|df| \leq 1$, this means 
\[ |df|(z) = 1 \quad \text{for all $z \in K$}.  \]
Then $f$ is a holomorphic isometric immersion from $K$ to $X$.
This contradicts the assumption that $X$ is $1/2$-WFL. 
\end{proof}

The above proposition states that there exists a constant $c(K) \lvertneqq 1$ satisfying the statement 
for \textit{each} unit square $K$. Since we have the Euclidean symmetry, we can adjust $c(K)$ so that
\[ c(K) = c(K') \quad \text{for all unit squares $K$ and $K'$}. \]
Hence we can define a universal constant $c \lvertneqq 1$ by setting
\[ c := c(K) \quad \text{for all unit squares $K$}.  \]
Then we can prove Theorem \ref{thm:upper-bound}.

\begin{proof}[\indent\sc Proof of Theorem \ref{thm:upper-bound}]
As it is noted first, we can suppose that $X$ is $1/2$-WFL by using a scale change. We have a constant $c \lvertneqq 1$ which satisfies the statement of Proposition \ref{prop:key} for all unit squares.

Let $f$ be an element of $\moduli (X)$ and $R$ be a positive number greater than $\sqrt{2}$.
We prove $\rho (f) \leq c$ by \textit{packing} unit squares in the disk $\bar{\Delta}(R) =
\{ z \in \affine | \, |z| \leq R \}$. (Here we promise that unit squares are closed.)

Since the diameter of a unit square is $\sqrt{2}$, we have the following fact: If a unit square $K$ has a intersection with the disk $\bar{\Delta}(R-\sqrt{2})$, $K$ is contained in $\bar{\Delta}(R)$.
Hence if we consider a tiling of the complex plane by unit squares, the disk $\bar{\Delta}(R-\sqrt{2})$
is covered by unit squares contained in $\bar{\Delta}(R)$. 
In other words, we have the following situation: 

There are unit squares $K_1, K_2, \dots, K_N$ contained in $\bar{\Delta}(R)$ such that 
\begin{gather*}
   \text{when $i \neq j$, $K_i$ and $K_j$ have common points at most on their boundaries,}\\
   \bar{\Delta}(R-\sqrt{2}) \subset \bigcup_{i = 1}^N K_i \subset \bar{\Delta}(R).   
\end{gather*}
$N = \mathrm{vol}(\bigcup K_i)$ satisfies 
\[ \pi (R - \sqrt{2})^2 \leq N \leq \pi R^2.  \] 
Therefore
\begin{equation*}
   \begin{split}
   \int_{\bar{\Delta}(R)} |df|^2 \, dxdy &= \sum_{i=1}^N \int_{K_i} |df|^2 \, dxdy
                                           + \int_{\bar{\Delta}(R) \setminus \cup Ki} |df|^2 \, dxdy,\\ 
                                         &\leq N \cdot c  +  (\pi R^2 - N) , \\
                                         &\leq \pi R^2 \cdot c  +  \{ \pi R^2 - \pi (R - \sqrt{2})^2 \}.
   \end{split}
\end{equation*}
Here we used Proposition \ref{prop:key} and $|df| \leq 1$. Dividing the above by $\pi R^2$, we have
\[ \frac{1}{\pi R^2} \int_{\bar{\Delta}(R)} |df|^2 \, dxdy
    \leq c + \{ 1 - (1 - \frac{\sqrt{2}}{R})^2 \}.  \]
Taking the superior limit, we get 
\[ \rho(f) = \limsup_{R \to \infty} \frac{1}{\pi R^2} \int_{\bar{\Delta}(R)} |df|^2 \, dxdy 
\leq c \lvertneqq 1. \]
Thus 
\[ \rho (X) \leq c \lvertneqq 1. \]
\end{proof}
\begin{remark}
The above argument does \textit{not} give an effective estimate of the constant $c$.
Hence we need another method if we want an explicit upper bound for $\rho(X)$. 
This is the theme of Section 5, and we prove Theorem \ref{thm:effective} there.
\end{remark}


\subsection{Lower bounds for holomorphic capacity}
Next we will establish lower bounds for holomorphic capacity by using holomorphic maps from 
elliptic curves to $X$. 

Let $\elliptic$ be a elliptic curve. Here $\Lambda$ is a lattice in the complex plane
$\affine$ with $\mathrm{rank}\, \Lambda = 2$. The elliptic curve $\elliptic$ has the Euclidean metric
induced by the universal covering $\affine$.

For a holomorphic map $f : \elliptic \to X$, we define its energy $E(f)$ by setting 
\begin{equation}
E(f) :=  \int_{\elliptic} f^* \omega = \int_{\elliptic} |df|^2 \, dxdy.  \label{def:energy}
\end{equation}
When $(X, \, \omega) = (\affine P^n , \, \omega_{FS})$, this is the degree of $f$ defined in
 (\ref{def:deg}). 
Let $\Tilde{f} : \affine \to X$ be the lift of $f$ and
 set $m := \sup_{z \in \elliptic} |df|(z)$.

\begin{proposition} \label{prop:elliptic}
\[ \lim_{R \to \infty} \frac{1}{\pi R^2} \int_{|z| \leq R} |d\Tilde{f}|^2 \, dxdy  = 
\frac{E(f)}{\mathrm{vol}( \elliptic ) }. \]
\end{proposition}
\begin{proof}
Let $\omega_1, \omega_2$ be a basis of $\Lambda$:
 $\Lambda = \mathbb{Z} \omega_1 \oplus \mathbb{Z} \omega_2$.
Let $L$ be a \textit{period parallelogram} of $\Lambda$, i.e.,
 $L = \{ z_0 + s \omega_1 + t \omega_2 \in \affine | \, 0\leq s, t \leq 1 \}$ for some $z_0 \in \affine$.
 We have $\mathrm{vol} (L) = \mathrm{vol} (\elliptic)$ and 
\[ \int_{L} |d\Tilde{f}|^2 \, dxdy  =  E(f) .  \]
We will use an argument similar to the proof of Theorem \ref{thm:upper-bound}. We prove the statement by \textit{packing} period parallelograms in the disk 
$\bar{\Delta}(R) = \{ z \in \affine | \, |z| \leq R \}$.

Let $l$ be the diameter of a period parallelogram. Then, if a period parallelogram $L$ intersects
with the disk $\bar{\Delta}(R - l)$, $L$ is contained in $\bar{\Delta}(R)$, ($R \geq l$).
Hence if we consider a tiling of the complex plane by period parallelograms, 
the disk $\bar{\Delta}(R - l)$ is covered by period parallelograms contained in $\bar{\Delta}(R)$.

In other words,
there are period parallelograms $L_1, L_2, \dots, L_N$ contained in $\bar{\Delta}(R)$ such that
\begin{gather*}
\text{when $i \neq j$, $L_i$ and $L_j$ have common points at most on their boundaries}, \\
\bar{\Delta}(R-l) \subset \bigcup_{i=1}^N  L_i  \subset \bar{\Delta}(R).
\end{gather*}
Then, $\mathrm{vol}(\bigcup L_i) = N \cdot \mathrm{vol}(\elliptic)$ satisfies 
\begin{equation} \label{estimate of period parallelograms}
\pi (R-l)^2 \leq \mathrm{vol}(\bigcup L_i) \leq \pi R^2. 
\end{equation} 
We have
\begin{equation*}
   \begin{split}
   \frac{1}{\pi R^2} \int_{\bar{\Delta}(R)} |d\Tilde{f}|^2 \, dxdy  &= 
          \frac{1}{\pi R^2} \sum_{i=1}^N \int_{L_i} |d\Tilde{f}|^2 \, dxdy   + 
          \frac{1}{\pi R^2} \int_{\bar{\Delta}(R) \setminus \cup L_i} |d\Tilde{f}|^2 \, dxdy  \\
   &= \frac{N}{\pi R^2} \cdot E(f) + 
       \frac{1}{\pi R^2} \int_{\bar{\Delta}(R) \setminus \cup L_i} |d\Tilde{f}|^2 \, dxdy.
   \end{split}
\end{equation*}
Hence
\begin{equation*}
  \begin{split}
   \left| \frac{1}{\pi R^2} \int_{\bar{\Delta}(R)} |d\Tilde{f}|^2 \, dxdy 
                                    - \frac{E(f)}{\mathrm{vol}(\elliptic)} \right|  &\leq
   \left| \frac{\mathrm{vol}(\cup L_i)}{\pi R^2} - 1 \right| \frac{E(f)}{\mathrm{vol}(\elliptic)}
   + \frac{m^2}{\pi R^2} \left| \pi R^2 - \mathrm{vol}(\cup L_i) \right|, \\
   &= \left( \frac{E(f)}{\mathrm{vol}(\elliptic)} + m^2 \right) \left| 1 - \frac{\mathrm{vol}(\cup L_i)}{\pi R^2} \right|. 
  \end{split}
\end{equation*}
Here $m = \norm{df}$. From the estimate of $\mathrm{vol}(\bigcup L_i)$ in (\ref{estimate of period parallelograms}), we have
\[ \left| 1 - \frac{\mathrm{vol}(\cup L_i)}{\pi R^2} \right|  \leq 1 - (1 - \frac{l}{R})^2  \to 0, 
\quad (R \to \infty). \]
Thus 
\[ \lim_{R \to \infty} \frac{1}{\pi R^2} \int_{|z| \leq R} |d\Tilde{f}|^2 \, dxdy  = 
\frac{E(f)}{\mathrm{vol}( \elliptic ) }. \]
\end{proof}

Then we can prove the main result of this section.
\begin{theorem} \label{thm:lower-bound}
If there exists a non-constant holomorphic map $f : \elliptic \to X$, we have the following estimate.
\[ \rho(X) \geq \frac{E(f)}{\norm{df}^2 \mathrm{vol}(\elliptic)}. \]
In particular,
\[ \rho (X) \gvertneqq 0. \]
\end{theorem}
\begin{proof}
Note that $m = \norm{df}$ and $E(f)$ are positive because $f$ is a non-constant map. 
From Lemma \ref{lem:scaling} and Proposition \ref{prop:elliptic},
\[  0 \lvertneqq \frac{E(f)}{\norm{df}^2 \mathrm{vol}(\elliptic)} = 
\lim_{R \to \infty} \frac{1}{m^2 \pi R^2} \int_{|z| \leq R} |d\Tilde{f}|^2 \, dxdy \leq
\rho(X). \]
\end{proof}

As a corollary of Theorem \ref{thm:lower-bound}, we get Theorem \ref{thm:gap}:
\begin{proof}[\indent\sc Proof of Theorem \ref{thm:gap}.]
If $f: \elliptic \to \affine P^n$ is a constant map, the statement is trivial.
When $f$ is a non-constant map, it follows from Theorem \ref{thm:lower-bound}.
(Here we have $E(f) = \mathrm{deg}(f)$.)
\end{proof}


\section{Proof of Theorem \ref{thm:non-trivial}}
We give the proof of Theorem \ref{thm:non-trivial} in this section.
Most of the proof have already been done in Section 2. The last piece which we need is the following.
\begin{proposition} \label{thm:Calabi}
The complex projective space $\affine P^n$ is WFL.
\end{proposition}
This proposition follows from a general result of E. Calabi concerning isometric imbeddings of complex manifolds [C, Theorem 8].
Here we will give a direct proof.

First, we recall the following easy fact:
\begin{lemma}\label{lem:lagrangian}
Let $f(z_1, z_2)$ be a holomorphic function in two variables defined on a connected open neighborhood of the origin in $\affine^2$. Set $g(z) := f(z, \bar{z})$. $g(z)$ is defined on a neighborhood of the origin in $\affine$. Then, if $g(z) \equiv 0$, we have $f(z_1, z_2) \equiv 0$.
\end{lemma}
\begin{proof}
Differentiating the equation $g(z) = f(z, \bar{z})$, we have
\begin{gather*}
\frac{\partial}{\partial z} g(z) :=
\frac{1}{2}\left( \frac{\partial}{\partial x} - \sqrt{-1} \frac{\partial}{\partial y} \right) g(z) =
\frac{\partial f}{\partial z_1} (z, \bar{z}), \\
\frac{\partial}{\partial \bar{z}} g(z) :=
\frac{1}{2}\left( \frac{\partial}{\partial x} + \sqrt{-1} \frac{\partial}{\partial y} \right) g(z) =
\frac{\partial f}{\partial z_2} (z, \bar{z}).
\end{gather*}
More generally we have
\[ \frac{\partial^{n+m}}{\partial z^n \partial \bar{z}^m} g(z) =
\frac{\partial^{n+m} f}{\partial z_1^n \partial z_2^m} (z, \bar{z}) 
\quad \text{for all $n, m \geq 0$.} \]

Hence, if $g(z) \equiv 0$, all partial derivatives of $f$ at $(z_1, z_2) = (0, 0)$ are zero.
This means that $f(z_1, z_2) \equiv 0$.
\end{proof}
\begin{proof}[\indent\sc Proof of Proposition \ref{thm:Calabi}]
Suppose that there exists a holomorphic isometric immersion $f: \Delta(r) \to \affine P^n$ 
for some positive number $r$. 
Because the complex projective space $\affine P^n$ is a homogeneous space, we can suppose that $f(0) = [1 : 0: \cdots: 0]$ without loss of 
generality. Then we can express $f$ in some neighborhood of the origin by
\[ f(z) = [1: f_1(z): f_2(z): \cdots : f_n(z)].  \]
Here $f_i(z)$ is a holomorphic function such that $f_i(0) = 0$, ($1 \leq i \leq n$).

The following argument is purely local. Hence we promise that all functions are defined on some
neighborhood of the origin in the complex plane $\affine$.

From the definition of the Fubini-Study metric (\ref{def:Fubini-Study}), we have
\[ |df|^2 (z) = \frac{1}{4 \pi} \Delta \log \left(1 + \sum_{i=1}^n |f_i(z)|^2 \right),  \quad 
(\Delta := \frac{\partial^2}{\partial x^2} + \frac{\partial^2}{\partial y^2} ). \]
Since $f$ is a holomorphic isometric immersion, $|df|^2 \equiv 1$. On the other hand, we have
\[ \frac{1}{4 \pi} \Delta (\pi |z|^2) \equiv 1. \]
Hence  
\[ \Delta \left\{ \log \left(1 + \sum |f_i(z)|^2 \right) - \pi |z|^2 \right\} \equiv 0. \] 
Because a harmonic function is locally the real part of a holomorphic function, 
we have a holomorphic function $g(z)$ such that $g(0) = 0$ and
\[ \log \left( 1 + \sum |f_i(z)|^2 \right) = \pi |z|^2 + g(z) + \overline{g(z)}. \]
Introducing the new holomorphic functions $\bar{f_i}(z) := \overline{f_i(\bar{z})}$ and
 $\bar{g}(z) := \overline{g(\bar{z})}$, we can express the above equation by 
\begin{equation*}
\log \left( 1 + \sum f_i(z) \bar{f_i}(\bar{z}) \right) = \pi z \bar{z} + g(z) + \bar{g}(\bar{z}).  
\end{equation*}
Applying Lemma \ref{lem:lagrangian} to this, we get
\begin{equation*}
\log \left( 1 + \sum f_i(z_1) \bar{f_i}(z_2) \right) = \pi z_1 z_2 + g(z_1) + \bar{g}(z_2).
\end{equation*}
Substituting $z_2 = 0$, we get $g(z_1) \equiv 0$ because $\bar{f_i}(0) = \bar{g}(0) =0$.
Thus the above equation becomes 
\begin{equation*}
\log \left( 1 + \sum f_i(z_1) \bar{f_i}(z_2) \right) = \pi z_1 z_2.
\end{equation*}
Hence we have 
\begin{equation*}
1 + \sum_{i=1}^n f_i(z_1) \bar{f_i}(z_2) = \exp (\pi z_1 z_2).
\end{equation*}
Applying $\partial^{\alpha + \beta} / \partial z_1^\alpha  \partial z_2^\beta$ 
at $(z_1, z_2) = (0, 0)$, we get
\begin{equation*}
\sum_{i=1}^n f_i^{(\alpha)}(0) \overline{f_i^{(\beta)}(0)} = \pi^{\alpha} \alpha ! \delta_{\alpha \beta} 
\quad \text{for $\alpha, \beta \geq 1$.}
\end{equation*} 
Here $f_i^{(\alpha)}(0)$ is the $\alpha$-th derivative of $f_i$ at the origin and 
$\delta_{\alpha \beta}$ is the Kronecker delta.
This means that \textit{an infinite number of} non-zero vectors 
$(f_1^{(\alpha)}(0),\, f_2^{(\alpha)}(0),\, \cdots,\, f_n^{(\alpha)}(0) )$ in $\affine^n$,
 $(\alpha \geq 1)$, are orthogonal to each other. It is impossible.
\footnote{This argument is essentially equivalent to the notion of ``resolvability of rank $n$'' in [C].} 
\end{proof}
\begin{proof}[\indent\sc Proof of Theorem \ref{thm:non-trivial}]
The complex projective space $\affine P^n$ satisfies 
all the conditions required in Theorem \ref{thm:upper-bound}. Hence we get
\[ \rho(\affine P^n) \lvertneqq 1. \]

Next we prove the lower bound. 
There are non-constant elliptic functions, for example, Weierstrass' elliptic function $\wp (z)$. 
Then we can apply Theorem \ref{thm:lower-bound}, and we get
\[ \rho(\affine P^1) \gvertneqq 0. \]
Since the natural inclusion $\affine P^1 \hookrightarrow \affine P^n$ is a holomorphic isometric 
imbedding, we conclude that
\[ \rho (\affine P^n) \geq \rho(\affine P^1) \gvertneqq 0. \]
\end{proof}


\section{Holomorphic capacity and theta functions}

We prove Theorem \ref{thm:limit} in this section. The proof is based on the lower bound given in Theorem \ref{thm:lower-bound}.
To apply it, we need appropriate holomorphic maps from elliptic curves to the complex projective spaces
$\affine P^n$. The classical theory of theta functions gives a good answer. 

The basic reference on theta functions is [M]. Actually, we will use only very basic facts on theta functions. All we need are given in Chapter 1, \S1, \S3 and \S4 in [M].

Let $\tau = s + t \sqrt{-1}$ be an element of the upper half plane, $(t = \mathrm{Im}\, \tau >0)$.
We define the theta function $\theta (z)$ by setting
\begin{equation*}
\theta (z) := \sum_{n \in \mathbb{Z}} \exp(\pi \sqrt{-1} n^2 \tau + 2\pi \sqrt{-1} n z)
\quad \text{for $z \in \affine$}.
\end{equation*}
Here we consider that $\tau$ is fixed, and we drop the dependence on $\tau$ in this notation for simplicity.

We define the theta function $\theta_{a, b} (z)$ with characteristics $a, b \in \mathbb{R}$ by setting
\begin{equation*}
\theta_{a,b}(z) := \exp (\pi \sqrt{-1} a^2 \tau + 2\pi \sqrt{-1} a (z+b)) \theta (z + a \tau +b).
\end{equation*}
We can construct projective imbeddings of elliptic curves by means of these $\theta_{a,b} (z)$. 
These projective imbeddings give good lower bounds for our problem.
We follow the arguments in [M, Chapter 1, \S4].

We define a lattice $\Lambda$ in $\affine$ by setting $\Lambda := \mathbb{Z} \oplus \mathbb{Z} \tau$.
For any integer $l \geq 2$, we set
\begin{gather*}
(\frac{1}{l} \mathbb{Z})^2 \cap [0, 1)^2 =  \{ (a_0,b_0), (a_1,b_1),\cdots, (a_{l^2-1}, b_{l^2-1}) \},\\
\theta_i (z) := \theta_{a_i, b_i} (z), \quad ( 0 \leq i \leq l^2-1 ),\\
\varphi_l : \affine/l \Lambda \to \affine P^{l^2-1}, \quad [z] \mapsto [\theta_0(z):\, \theta_1(z):\, \cdots:\, \theta_{l^2-1}(z)].
\end{gather*}
The following fact is proved in [M, Chapter 1, \S4]. (The statement on the degree of $\varphi_l$ follows from Lemma 4.1 there.)
\begin{fact}
$\varphi_l$ is a well-defined holomorphic imbedding and $\mathrm{deg} (\varphi_l) = l^2$.
\end{fact}
Since $\mathrm{vol} (\affine / l \Lambda) = t\, l^2$, ($t = \mathrm{Im} \, \tau$), 
Theorem \ref{thm:lower-bound} for $\varphi_l$ gives the following lower bound:
\begin{equation} \label{estimate by theta}
\rho( \affine P^{l^2-1}) \geq \frac{1}{t  \norm{d\varphi_l}^2}.
\end{equation}
Therefore we need the estimate of $\norm{d\varphi_l}$. The following lemma is the basis of our argument.
\begin{lemma}\label{lem:limit}
\[ \int_0^1 \int_0^1 |\theta_{a,b}(z)|^2\, dadb =
 \sqrt{\frac{1}{2 t}} \exp \left( \frac{2\pi y^2}{t}\right). \]
 Hence we have
\[ \frac{1}{4\pi}\Delta \log \int_0^1 \int_0^1 |\theta_{a,b}(z)|^2\, dadb  \equiv  \frac{1}{t},
\quad (\Delta = \frac{\partial^2}{\partial x^2} + \frac{\partial^2}{\partial y^2}). \]
\end{lemma}
\begin{proof}
The proof is just a calculation.
\begin{gather*}
 |\theta_{a,b} (z)| = \exp (-\pi a^2 t - 2\pi a y) |\theta (z + a \tau +b)|, \\
 \theta (z + a\tau + b) = 
\sum_n \exp (\pi \sqrt{-1} n^2 \tau + 2\pi \sqrt{-1} n(z + a\tau )) \exp (2\pi \sqrt{-1} n b).  
\end{gather*}
From Parseval's equality 
\begin{equation*}
 \begin{split}
 \int_0^1 |\theta (z + a\tau + b)|^2 db  &= 
 \sum_n |\exp (\pi \sqrt{-1} n^2\tau + 2\pi \sqrt{-1} n (z + a\tau )) |^2 , \\
 &= \sum_n \exp (-2\pi n^2 t - 4\pi n(y+ at) ).
 \end{split}  
\end{equation*}
Hence 
\begin{equation*}
 \begin{split}
 \int_0^1 |\theta_{a,b} (z)|^2 \, db  &=
 \exp (-2\pi a^2 t -4\pi a y) \sum_n \exp (-2\pi n^2 t - 4\pi n(y+ at) ), \\
 &= \sum_n \exp (-2\pi t (a+n)^2  - 4\pi y (a+n)).
 \end{split}
\end{equation*}
Thus we get
\begin{equation*}
 \begin{split}
 \int_0^1 da \int_0^1 |\theta_{a,b} (z)|^2 db 
 &= \sum_n \int_0^1 \exp (-2\pi t (a+n)^2  - 4\pi y (a+n))\, da , \\
 &= \sum_n \int_n^{n+1}  \exp (-2\pi t a^2  -4\pi y a) \, da , \\
 &= \int_{-\infty}^{+\infty}  \exp (-2\pi t a^2 -4\pi y a) \, da , \\
 &= \sqrt{\frac{1}{2 t}} \exp \left( \frac{2\pi y^2}{t}\right) .
 \end{split}
\end{equation*}
\end{proof} 
We need one more lemma.
\begin{lemma} \label{lem:equivariance}
For $\alpha, \beta \in \mathbb{Z}$, we have 
\[ |d\varphi_l | (z + \alpha \tau + \beta) = |d\varphi_l |(z). \]
In other words, $|d\varphi_l|(z)$ is invariant under the following $\mathbb{Z}^2$-action on 
$\affine / l \Lambda$.
\[ \mathbb{Z}^2 \curvearrowright \affine/ l \Lambda, 
\quad ((\alpha, \beta), z) \mapsto z + \alpha \tau + \beta. \]
\end{lemma}
\begin{proof}
\[ \theta_{a,b} (z + \alpha \tau + \beta) =
\exp (-\pi \sqrt{-1} \alpha^2 \tau - 2\pi \sqrt{-1} \alpha z) 
\exp (-2\pi \sqrt{-1} \alpha b)\, \theta_{a+\alpha, b+ \beta} (z). \]
Since $\alpha$ and $\beta$ are integers, we have 
\[ \theta_{a + \alpha, b + \beta} (z) = \exp (2\pi \sqrt{-1} a \beta )\, \theta_{a, b} (z). \]
Set $c_i := \exp (-2\pi \alpha b_i + 2\pi \sqrt{-1} a_i \beta) \in U(1)$, ($0 \leq i \leq l^2 -1$).
Then 
\[ \theta_i (z + \alpha \tau + \beta ) =
\exp (-\pi \sqrt{-1} \alpha^2 \tau - 2\pi \sqrt{-1} \alpha z)\, c_i \, \theta_i (z). \]
It results that 
\[ \varphi_l (z+ \alpha \tau + \beta ) =
 [ c_0\, \theta_0 (z): \, c_1\, \theta_1 (z): \cdots : c_{l^2-1}\, \theta_{l^2-1}(z)]. \]
(This is the ``equivariance'' described in [M, Chapter 1, \S4]. )

Hence we get
\begin{equation*}
 \begin{split}
 |d\varphi_l |^2 (z + \alpha \tau + \beta ) &=
  \frac{1}{4\pi} \Delta \log \sum_i |c_i\, \theta_i (z)|^2,  \\
 &= \frac{1}{4\pi} \Delta \log \sum_i |\theta_i (z)|^2, \\
 &= |d\varphi_l |^2 (z).
 \end{split}
\end{equation*} 
\end{proof}

\begin{remark}
Lemma \ref{lem:equivariance} is still true for $\alpha , \beta \in \frac{1}{l} \mathbb{Z}$,
but we don't need it.
\end{remark}

The following proposition gives a sufficient estimate for the proof of Theorem \ref{thm:limit}. 
This is a special case of the results of G. Tian [T].
\begin{proposition} \label{prop:Tian} 
\[ \lim_{l \to \infty} \norm{|d\varphi_l |^2 - \frac{1}{t}} = 0. \]
Hence we have
\[ \lim_{l \to \infty} \norm{d\varphi_l}^2 = \frac{1}{t}. \]
\end{proposition}
\begin{proof}
From the definition of the Fubini-Study metric (\ref{def:Fubini-Study}), 
\[ |d\varphi_l|^2 (z) = \frac{1}{4\pi} \Delta \log \sum_i |\theta_i (z)|^2
 = \frac{1}{4\pi} \Delta \log \left( \frac{1}{l^2}\sum_i |\theta_i(z)|^2 \right) . \]
Since $[0,1)^2 = \bigsqcup_{i} [a_i, a_i + 1/l) \times [b_i, b_i + 1/l)$ is a division into small squares, 
the definition of the Riemann integral gives the following \textit{point-wise} convergence.
\[ \lim_{l \to \infty}  \frac{1}{l^2}\sum_i |\theta_i(z)|^2 =
  \int_0^1 \int_0^1 |\theta_{a,b} (z)|^2 \, dadb 
  = \sqrt{\frac{1}{2 t}} \exp \left( \frac{2\pi y^2}{t}\right)  
\quad \text{for any $z \in \affine$}. \]
Actually we can say more. Set $K := \{ x + y \tau \in \affine |\, 0 \leq x, y \leq 1 \}$.
Since $K$ is compact, it is easy to see that 
\[ \lim_{l \to \infty} \Dnorm{\frac{1}{l^2}\sum_i |\theta_i(z)|^2 - \int_0^1 \int_0^1 |\theta_{a,b} (z)|^2 \, dadb}{k}{K} = 0  \quad \text{for all $k \geq 0$}. \]
Here $\Dnorm{\, \cdot \,}{k}{K}$ is the $\mathcal{C}^k$-norm for functions defined over $K$.
Therefore we have
\[ \lim_{l \to \infty} \Dnorm{\frac{1}{4\pi} \Delta \log \left( \frac{1}{l^2}\sum_i |\theta_i(z)|^2 \right) -  \frac{1}{4\pi}\Delta \log \int_0^1 \int_0^1 |\theta_{a,b}(z)|^2\, dadb}{0}{K} = 0. \]  
Hence 
\[ \lim_{l \to \infty} \Dnorm{ |d\varphi_l |^2 - \frac{1}{t} }{0}{K} = 0. \]  
Here we consider $K$ as a subspace of the elliptic curve $\affine / l \Lambda$ 
through the natural projection $\affine \to \affine / l \Lambda$. 

Since $K$ is a fundamental domain for the $\mathbb{Z}^2$-symmetry described in Lemma \ref{lem:equivariance}, we get the conclusion:
\[ \lim_{l \to \infty} \Dnorm{ |d\varphi_l |^2 - \frac{1}{t} }{0}{\affine / l \Lambda} = 0. \]
\end{proof}

\begin{proof}[\indent\sc Proof of Theorem \ref{thm:limit}]
From the inequality (\ref{estimate by theta}) and Proposition \ref{prop:Tian}, we have
\[ \liminf_{l \to \infty} \rho (\affine P^{l^2-1}) \geq 
\lim_{l \to \infty} \frac{1}{t  \norm{d\varphi_l}^2} = 1. \]
Hence $\lim_{l \to \infty} \rho (\affine P^{l^2-1}) = 1$. 
On the other hand, we have
\[ \rho(\affine P^1) \leq \rho(\affine P^2) \leq \rho(\affine P^3) \leq \cdots
\leq \rho(\affine P^n) \leq \rho(\affine P^{n+1}) \leq \cdots \lvertneqq 1. \]
Thus  
\[ \lim_{n \to \infty} \rho (\affine P^n) = 1. \]
\end{proof}


\section{Explicit upper bounds for $\rho(\affine P^1)$}
In this section we prove Theorem \ref{thm:effective} which gives an explicit upper bound on the holomorphic capacity of the complex projective line.

The complex projective line is the Riemann sphere: $\affine P^1 = \affine \cup \{ \infty \}$. Hence a holomorphic map
from $\affine$ to $\affine P^1$ is a meromorphic function in the complex plane.

To simplify the calculations, we will use the following rescaled Fubini-Study metric $ds^2$
as the metric on $\affine P^1$ in this section. (A holomorphic capacity is invariant under a scale change of a metric, (Proposition \ref{prop:scaling}).)
\begin{gather}
ds^2 := \frac{dw d\bar{w}}{(1+ |w|^2)^2}, \label{def:rescaled Fubini-Study} \\
\text{the fundamental 2-form of $ds^2$} = \frac{\sqrt{-1}}{2}\frac{dw \wedge d\bar{w}}{(1+|w|^2)^2}
\quad (= \pi \omega_{FS}). \quad \text{(cf. (\ref{def:Fubini-Study}).)} \notag
\end{gather}  
Here $w$ is the natural coordinate on $\affine$ of $\affine P^1 = \affine \cup \{ \infty \}$.

Then, for a meromorphic function $f(z)$ in the complex plane, the norm of the differential $df$ is expressed by 
\begin{equation*} 
|df|(z) = \frac{|f'(z)|}{1+ |f(z)|^2}.
\end{equation*}


\subsection{Preliminary estimates} \label{preliminary estimates}
In this section we prepare various estimates for the proof of Theorem \ref{thm:effective}.
Here we do \textit{not} pursue precise estimates. Actually we will use many loose estimates for simplicity of the calculations. Our purpose is to show the fact that we can get an explicit upper bound on a holomorphic capacity.

First we compute the distance on $\affine P^1$ defined by the metric (\ref{def:rescaled Fubini-Study}).
\begin{lemma} \label{lem:distance of rescaled Fubini-Study}
Let $w$ be a point of $\affine P^1$, then we have 
\[ d_{\affine P^1} (0, w) = \arctan |w|.  \]
Here $d_{\affine P^1} (\cdot, \cdot)$ is the distance on $\affine P^1$ defined by $(\ref{def:rescaled Fubini-Study})$ and $\arctan (\cdot)$ is the branch of the inverse function of $\tan (\cdot)$ satisfying $\arctan 0 = 0$.
\end{lemma}
\begin{proof}
It is easy to see that $c(t) := w t,\, (0 \leq t \leq 1)$, is the minimum geodesic from $0$ to $w$.
Hence we get 
\[ d_{\affine P^1} (0, w) = \int_0^1 \frac{|c'(t)|}{1 + |c(t)|^2} dt 
= \int_0^1 \frac{|w|}{1 + |w|^2 t^2} dt = \arctan |w|. \]
\end{proof}

We set $r:= |z|$ for $z$ in the complex plane and put
 \[ \varepsilon := 10^{-100}\quad  \text{and}\quad  r_0 := 10^{-10} . \]

Let $f(z)$ be a meromorphic function in the complex plane which satisfies 
\begin{equation}\label{rescaled norm of df}
 |df|(z) = \frac{|f'(z)|}{1 + |f(z)|^2} \leq 1 \quad \text{for all $z\in \affine$}. 
\end{equation}
In addition, we suppose that the following conditions are also satisfied:
\begin{equation}\label{pre-estimate}
f(0) = 0 \quad \text{and}\quad  1 - \varepsilon \leq |df|(0) \leq 1.
\end{equation}

The main purpose of this section is to get a good estimate of $|df|^2(z)$ in the small disk $r  \leq r_0$ under the conditions (\ref{rescaled norm of df}) and (\ref{pre-estimate}).

The following lemma is the basis of our argument.
\begin{lemma} \label{lem:apriori estimate}
\[ |f(z)|^2 \leq \tan ^2 r \leq r^2 + 2 r^4 \leq 2 r^2, \quad (r \leq r_0). \]
Let $f(z) = a_1 z + a_2 z^2 + a_3 z^3 + \cdots$ be the Taylor expansion centered at the origin. Then we have
\[ 1 - \varepsilon \leq |a_1| \leq 1 \quad \text{and} \quad |a_n| \leq (4/\pi )^n < 2^n. \]
\end{lemma}
\begin{proof}
Since $|df| \leq 1$ and $f(0) = 0$, we have 
\[ \arctan |f(z)| = d_{\affine P^1} (f(0),f(z)) \leq r. \] 
Hence 
\begin{equation}
 |f(z)| \leq \tan r, \quad (r < \pi /2). \label{apriori estimate}
\end{equation}
Here we have
\[ \tan ^2 r = \frac{\sin ^2 r}{1 - \sin ^2 r} = \sin ^2 r +  \frac{\sin ^4 r}{1 - \sin ^2 r}
\leq r^2 + 2 r^4 \leq 2 r^2, \quad (r \leq r_0). \]
Hence we get the above first statement. (Of course, this is a \textit{very} loose estimate.)

Next we will estimate the coefficients of the Taylor expansion. 
Since $|df|(0) = |f'(0)| = |a_1|$, we have 
\[ 1 - \varepsilon \leq |a_1| \leq 1. \]
Using (\ref{apriori estimate}) at $r = \pi /4$, we have 
\[ |f(z)| \leq 1, \quad (r = \pi /4). \]
Thus we get
\[ |a_n| = \left| \frac{1}{2\pi \sqrt{-1}} \int_{|z| = \pi /4} \frac{f(z)}{z^{n+1}}\, dz \right| 
        \leq (4 / \pi)^n < 2^n. \]
\end{proof}

In the following estimates, we always assume $r \leq r_0$. The important term in the estimate of $|df|^2$ is the second order term and the higher order terms will be loosely estimated.

First we will estimate the denominator of $|df|^2(z) = |f'(z)|^2/(1+ |f(z)|^2)^2$: 
\begin{lemma} \label{lem:estimate of denominator}
\[ \frac{1}{(1+ |f(z)|^2 \,)^2} \leq 1 - 2\, r^2 + 30\,  r^3 + \varepsilon , \quad ( r \leq r_0). \]
\end{lemma}
\begin{proof}
Using Lemma \ref{lem:apriori estimate}, if $r \leq r_0$,  we have
\begin{equation*}
 \begin{split}
 \frac{1}{(1+ |f(z)|^2)^2} &= 1 - 2\,|f(z)|^2 + \frac{3\, |f(z)|^4 + 2\, |f(z)|^6}{(1+ |f(z)|^2)^2 }, \\
                           &\leq 1 - 2\,|f(z)|^2 + 3\cdot 2^2 \cdot r^4 + 2\cdot 2^3 \cdot r^6 , \\
                           &= 1 - 2\, |f(z)|^2 + ( 3\cdot 2^2 + 2^4 \cdot r^2 ) r^4, \\
                           &\leq  1 - 2\, |f(z)|^2 + 13\, r^4. 
 \end{split}
\end{equation*}
The term $-2\,|f(z)|^2$ can be estimated as follows.
\[ |f(z)| \geq |a_1 z| - |a_2 z^2 + a_3 z^3 + \cdots |. \]
And we have
\[ |a_2 z^2 + a_3 z^3 + \cdots | \leq 2^2\cdot r^2 + 2^3 \cdot r^3 + \cdots
 = \frac{4 r^2}{1 - 2 r} \leq 5\, r^2 . \]
Hence 
\[ |f(z)| \geq |a_1|\, r - 5\, r^2 \geq (1 - \varepsilon )\, r - 5\, r^2 \geq 0. \]
Then 
\[ |f(z)|^2 \geq (1- \varepsilon )^2\, r^2 - 10\, (1- \varepsilon )\, r^3 + 25 \, r^4 
            \geq (1 - 2 \, \varepsilon )\, r^2 - 10 \, r^3 . \]
Thus 
\[ -2 \, |f(z)|^2 \leq -2\, (1 - 2 \, \varepsilon )\, r^2 + 20 \, r^3 . \]
Hence we get the conclusion:
\begin{equation*}
 \begin{split}
 \frac{1}{(1+ |f(z)|^2 \, )^2} & \leq 1 - 2\, (1 - 2\, \varepsilon )\, r^2 + 20\, r^3 + 13 \, r^4 \\
                            & = 1 - 2\, r^2 + (20 + 13\, r)\, r^3 + 4 \, \varepsilon \, r^2  \\
                            & \leq 1 - 2\, r^2 + 30\, r^3 + \varepsilon.
 \end{split}
\end{equation*}
\end{proof}

Next we will estimate the numerator $|f'(z)|^2$. This case is more complicated and we need some 
preparations. 
\begin{lemma}\label{lem:estimate of a_2}
\[2 \, |a_2| \leq 30 \, \sqrt{\varepsilon}. \]
\end{lemma}
\begin{proof}
\[ f'(z) = a_1 + 2 a_2  z + 3 a_3  z^2 + 4 a_4  z^3 + \cdots . \]
Using $|df| \leq 1$, if $r\leq r_0$, we have 
\[ |f'(z)| \leq 1 + |f(z)|^2 \leq 1 + 2 \, r^2.  \]
Hence 
\begin{equation*}
 \begin{split}
 |a_1 + 2\, a_2 \, z| &\leq 1 + 2 \, r^2 + |3 a_3 z^2 + 4 a_4 z^3 + 5 a_5 z^4 + \cdots |, \\
 &\leq 1 + 2 \, r^2 + (3\cdot 2^3\cdot r^2 + 4\cdot 2^4\cdot r^3 + 5 \cdot 2^5 \cdot r^4 + \cdots ),  \\
 &= 1 + 2 \, r^2 + \frac{24\, r^2 - 32\, r^3}{(1 - 2r )^2}, \\
 &\leq 1 + r^2  \left( 2 + \frac{24}{(1 - 2r )^2} \right) , \\
 &\leq 1 + 27 \, r^2.
 \end{split}
\end{equation*} 

Let $\theta_n \in \mathbb{R}/ 2\pi \mathbb{Z}$ be the argument of $a_n$: 
$a_n = |a_n| e^{\sqrt{-1} \theta_n}$. If $a_n = 0$, we promise $\theta_n := 0$.

Putting $z = \sqrt{\varepsilon}\, e^{\sqrt{-1} (\theta_1 - \theta_2 )} $ in the above, we get
\[ |a_1| + 2\, |a_2|\, \sqrt{\varepsilon} \leq 1 + 27 \, \varepsilon .\] 
Since $|a_1| \geq 1 - \varepsilon$,
\[ 2\, |a_2| \leq \frac{ \varepsilon + 27 \, \varepsilon}{\sqrt{\varepsilon}} 
  \leq 30 \, \sqrt{\varepsilon} .\]
\end{proof}

Set $E(r) := 4\, |a_4| \, r^3 + 5 \, |a_5|\, r^4 + 6 \, |a_6|\, r^5 + \cdots$. 
(``E'' is the initial letter of ``error term''.)
\begin{lemma}\label{lem:estimate of error term}
\[ E(r) \leq 65 \, r^3 \leq 10^{-20}, \quad (r \leq r_0). \]
\end{lemma}
\begin{proof}
Using Lemma \ref{lem:apriori estimate}, if $r \leq r_0$, we have
\begin{equation*}
 \begin{split}
 E(r) &\leq 4\cdot 2^4 \cdot r^3 + 5\cdot 2^5 \cdot r^4 + 6\cdot 2^6 \cdot r^5 + \cdots , \\
      &= \frac{64 \, r^3 - 96 \, r^4}{(1 - 2 \, r)^2} 
      \leq \frac{64\, r^3 }{(1 - 2 \, r)^2} \leq  65 \, r^3 . 
 \end{split}
\end{equation*}
Since $r_0 = 10^{-10}$, we have $65\, r^3 < 10^{-20}$.
\end{proof}
\begin{lemma}\label{lem:estimate of a_3}
Set $\delta := 10^{-5}$. We have
\[ 3\, |a_3|\, r^2 \leq r^2 + 100\, r^3 + \delta \sqrt{\varepsilon}, \quad (r \leq r_0). \]
\end{lemma}
\begin{proof}
From Lemma \ref{lem:apriori estimate},
\[  |a_1 + 2 a_2  z + 3 a_3  z^2 + 4 a_4  z^3 + \cdots |
 = |f'(z)| \leq 1 + |f(z)|^2 \leq 1 + r^2 + 2\, r^4. \]
Using Lemma \ref{lem:estimate of a_2} and Lemma \ref{lem:estimate of error term}, we have
\begin{equation*}
 \begin{split}
 |a_1 + 3\, a_3 \, z^2| &\leq 1+r^2+ 2\, r^4 + |2\, a_2 \, z + 4\, a_4 \, z^3 + 5\, a_5 \, z^4 + \cdots |,\\
                        &\leq 1+ r^2 + 2\, r^4 + 30 \sqrt{\varepsilon}\, r + E(r) ,\\
                        &\leq 1 + r^2 + (2\, r + 65)\, r^3 + 30 \sqrt{\varepsilon}\, r ,\\
                        &\leq 1 + r^2 + 100\, r^3 + 30  \sqrt{\varepsilon} \, r .
 \end{split}
\end{equation*}

Let $\theta_n$ be the argument of $a_n$ as in the proof of Lemma \ref{lem:estimate of a_2}.
Putting $z = r \, e^{\frac{\sqrt{-1}}{2}(\theta_1 - \theta_3)}$ in the above, we get
\[ |a_1| + 3 \, |a_3|\, r^2 \leq 1 + r^2 + 100\, r^3 + 30  \sqrt{\varepsilon} \, r. \]
Since $|a_1| \geq 1 - \varepsilon$, we have
\[ 3 \, |a_3|\, r^2 \leq \varepsilon + r^2 + 100\, r^3 + 30 \sqrt{\varepsilon} \, r 
                    \leq r^2 + 100 \, r^3 + \delta \, \sqrt{\varepsilon}.  \]
\end{proof}

Set 
\[\cos^+x := \max \, (0,\, \cos x) \quad \text{for all $x \in \mathbb{R}$} . \]

\begin{lemma} \label{lem:estimate of numerator}
\[ |f'(z)|^2 \leq 1 + 2\, r^2 \, \cos^+(2\, \theta - \theta_1 + \theta_3) 
   + 500\, r^3 + \frac{1}{2} \sqrt{\varepsilon} , \quad (r \leq r_0).  \]
Here $\theta$ is the argument of $z$ and $\theta_n$ is the argument of $a_n$.
\end{lemma}
\begin{proof}
Since $f'(z) = a_1 + 2a_2 z + 3a_3 z^2 + \cdots$, we have
\begin{equation*}
 \begin{split}
 |f'(z)| &\leq |a_1 + 3\, a_3\, z^2 | + 
               |2\, a_2\, z|  + | 4\, a_4\, z^3 + 5\, a_5\, z^4 + 6\, a_6\, z^5 + \cdots |, \\
         &\leq |a_1 + 3\, a_3 \, z^2 | + 30\, \sqrt{\varepsilon}\, r + E(r).  
 \end{split}
\end{equation*}
Hence
\begin{equation*}
 \begin{split}
 |f'(z)|^2 \leq & |a_1 + 3\, a_3\, z^2 |^2 + E(r)^2 + 2\, |a_1 + 3\, a_3\, z^2 |\cdot E(r) \\
                & + 900\, \varepsilon \, r^2 + 2\, |a_1 + 3\, a_3 \, z^2| \cdot 30 \sqrt{\varepsilon}\, r
                + 60  \sqrt{\varepsilon} \, r \cdot E(r).
 \end{split}
\end{equation*}
Since $\, 2\, |a_1 + 3\, a_3 \, z^2| \leq 2 |a_1| + 6\, |a_3|\, r^2 \leq 2 + 6\cdot 2^3\cdot r^2 \leq 3$, 
\begin{equation*}
 \begin{split}
 900\, \varepsilon \, r^2 + 2\, |a_1 + 3\, a_3 \, z^2| \cdot 30 \sqrt{\varepsilon}\, r        
  + 60  \sqrt{\varepsilon} \, r \cdot E(r) 
  &\leq (900\, \sqrt{\varepsilon}\, r^2 + 90 \, r + 60 \, r\, E(r))\sqrt{\varepsilon} , \\
  &\leq \delta \sqrt{\varepsilon}, \quad ( \delta = 10^{-5} ).
 \end{split}
\end{equation*}
Here we used $r \leq r_0 = 10^{-10}$ and Lemma \ref{lem:estimate of error term}: $E(r) \leq 10^{-20}$.

Thus we get
\begin{equation*}
 \begin{split}
 |f'(z)|^2 &\leq |a_1 + 3\, a_3\, z^2 |^2 + E(r)^2 + 3\, E(r) + \delta \, \sqrt{\varepsilon} , \\
           &\leq |a_1 + 3\, a_3\, z^2 |^2 + 4\, E(r) + \delta \, \sqrt{\varepsilon} .
 \end{split}
\end{equation*}
Using Lemma \ref{lem:estimate of a_3}, we have
\begin{equation*}
 \begin{split}
 |a_1 + 3\, a_3 \, z^2 |^2 &=
 |a_1|^2 + 2\cdot 3\cdot |a_1|\cdot |a_3|\cdot r^2 \cos (2\theta - \theta_1 + \theta_3 ) + 9\, |a_3|^2 r^4,\\
 &\leq 1 + 2\cdot 3\, |a_3|\, r^2 \cos^+ (2\theta - \theta_1 + \theta_3 ) + 9\cdot 2^6 \, r^4 ,\\
 &\leq 1 + 2\, (r^2 + 100\, r^3 + \delta \, \sqrt{\varepsilon}) \cos^+ (2\theta - \theta_1 + \theta_3 )
 + 576\, r^4 , \\
 &\leq 1 + 2\, r^2 \cos^+ (2\theta - \theta_1 + \theta_3 ) + 200\, r^3 + 2\delta \cdot \sqrt{\varepsilon}
 + r^3 , \\
 & = 1 + 2\, r^2 \cos^+ (2\theta - \theta_1 + \theta_3 ) + 201\, r^3 + 2\delta \cdot \sqrt{\varepsilon}.
 \end{split}
\end{equation*}
Thus 
\begin{equation*}
 \begin{split}
 |f'(z)|^2 &\leq 1 +  2\, r^2 \cos^+ (2\theta - \theta_1 + \theta_3 ) + (201\, r^3 + 4\, E(r)) + 
 3\delta \cdot \sqrt{\varepsilon}, \\
 &\leq 1 +  2\, r^2 \cos^+ (2\theta - \theta_1 + \theta_3 ) + 500\, r^3 + \frac{1}{2}\sqrt{\varepsilon}.
 \end{split}
\end{equation*}
\end{proof}

Then we can estimate $|df|^2(z)$ in the small disk $r \leq r_0$.
\begin{proposition}\label{prop:estimate of the norm}
\[ |df|^2(z) \leq 1 - 2\, r^2\, (1 - \cos^+ (2\theta - \theta_1 + \theta_3 )) 
   + 600\, r^3 + \sqrt{\varepsilon} , \quad (r \leq r_0). \]
\end{proposition}
\begin{proof}
To simplify the descriptions, we set $\varphi := 2\theta - \theta_1 + \theta_3$.

From Lemma \ref{lem:estimate of numerator}, 
\begin{equation*}
 \begin{split}
|df|^2(z) =  \frac{|f'(z)|^2}{(1+|f(z)|^2)^2} &\leq \frac{1 + 2\, r^2\, \cos^+ \varphi  
   + 500\, r^3 + \frac{1}{2} \sqrt{\varepsilon}}{(1+|f(z)|^2)^2}, \\
 &\leq \frac{1 + 2\, r^2\, \cos^+ \varphi }{(1+|f(z)|^2)^2}
 + 500\, r^3 +  \frac{1}{2} \sqrt{\varepsilon}. 
 \end{split}
\end{equation*}
From Lemma \ref{lem:estimate of denominator},
\begin{equation*}
 \begin{split}
 \frac{1 + 2\, r^2\, \cos^+ \varphi }{(1+|f(z)|^2)^2} \leq & 
\, (1 + 2\, r^2\, \cos^+ \varphi )(1 - 2\, r^2 + 30 \, r^3 + \varepsilon ) , \\
 =& \, 1 - 2\, r^2 \, (1 - \cos^+ \varphi ) 
+ (30\, r^3 - 4\, r^4 \cos^+ \varphi + 60 \, r^5 \cos^+ \varphi ) \\
& + \varepsilon + 2\, \varepsilon \, r^2\, \cos^+ \varphi , \\
\leq & \,  1 - 2\, r^2 \, (1 - \cos^+ \varphi ) + 31\, r^3 + 2\, \varepsilon .
  \end{split}
\end{equation*}
Thus
\begin{equation*}
 \begin{split}
 |df|^2(z) &\leq 1 - 2\, r^2 \, (1 - \cos^+ \varphi ) + 531\, r^3 
+ (2\, \sqrt{\varepsilon} + \frac{1}{2} ) \sqrt{\varepsilon} , \\
 &\leq 1 - 2\, r^2 \, (1 - \cos^+ \varphi ) + 600\, r^3 + \sqrt{\varepsilon} .
 \end{split}
\end{equation*}
\end{proof}
The following proposition is the conclusion of Section \ref{preliminary estimates}.
\begin{proposition} \label{prop:estimate on sector}
Let $D$ be a circular sector of radius $r_0$ and angle $\pi /2$ centered at the origin in the complex plane, i.e.,
\[ D = \{ r\, e^{\sqrt{-1} \theta} \in \affine \, |\,
 0 \leq r \leq r_0,\> \alpha \leq \theta \leq \alpha + \pi /2 \, \}
 \quad \text{for some $\alpha \in \mathbb{R}$}. \]
Then 
\[ \frac{1}{\mathrm{vol}(D)} \int_D |df|^2 \, dxdy \leq 1 - \frac{1}{4}\, r_0^2 + \sqrt{\varepsilon}
  < 1 - 10^{-30}. \]
\end{proposition}
\begin{proof}
From Proposition \ref{prop:estimate of the norm},
\[ \frac{1}{\mathrm{vol}(D)} \int_D |df|^2 \, dxdy \leq 1 + \sqrt{\varepsilon} 
- \frac{2}{\mathrm{vol}(D)} \int_D r^2 \, (1 - \cos^+ (2\theta - \theta_1 + \theta_3 ) )
+ \frac{600}{\mathrm{vol}(D)} \int_D r^3 .\]
Since $\mathrm{vol}(D) = \pi r_0^2 /4$,
\[ \frac{1}{\mathrm{vol}(D)} \int_D r^3 = 
\frac{4}{\pi r_0^2} \int_0^{r_0} r^4\, dr \, \int_{\alpha}^{\alpha + \pi /2} d\theta
= \frac{2}{5}\, r_0^3 . \]
\begin{equation*}
 \begin{split}
 \frac{1}{\mathrm{vol}(D)} \int_D r^2 \, (1 - \cos^+ (2\theta - \theta_1 + \theta_3 ) )
 &= \frac{4}{\pi r_0^2}\int_0^{r_0} r^3 dr 
 \int_{\alpha}^{\alpha + \pi /2} (1 - \cos^+ (2\theta - \theta_1 + \theta_3 ) )  d\theta ,\\
 &= \frac{r_0^2}{\pi} \int_{\alpha}^{\alpha + \pi /2} 
(1 - \cos^+ (2\theta - \theta_1 + \theta_3 ) ) \, d\theta .
 \end{split}
\end{equation*}
Setting $\varphi := 2\theta -\theta_1 + \theta_3$ and $\beta := 2\alpha - \theta_1 + \theta_3$, we have
\begin{equation*}
 \begin{split}
 \frac{r_0^2}{\pi} \int_{\alpha}^{\alpha + \pi /2} 
(1 - \cos^+ (2\theta - \theta_1 + \theta_3 ) ) \, d\theta 
&= \frac{r_0^2}{2 \pi}\int_{\beta}^{\beta + \pi} (1 - \cos^+ \varphi )\, d\varphi , \\
&\geq \frac{r_0^2}{2 \pi}\int_{-\pi /2}^{\pi /2} (1 - \cos^+ \varphi )\, d\varphi , \\
&= r_0^2 \, (\frac{1}{2} - \frac{1}{\pi}). 
 \end{split}
\end{equation*}
Thus 
\begin{equation*}
 \begin{split}
 \frac{1}{\mathrm{vol}(D)} \int_D |df|^2 \, dxdy &\leq 
 1 + \sqrt{\varepsilon} - (1- \frac{2}{\pi})\, r_0^2 + 240\, r_0^3, \\
 &\leq 1 + \sqrt{\varepsilon} - \frac{1}{3}\, r_0^2 + 240\, r_0^3, \\
 &\leq 1 -\frac{1}{4}\, r_0^2 + \sqrt{\varepsilon}.
 \end{split}
\end{equation*}
Since $r_0 = 10^{-10}$ and $\varepsilon = 10^{-100}$, 
\[ 1 -\frac{1}{4}\, r_0^2 + \sqrt{\varepsilon} = 1 - \frac{10^{-20}}{4} + 10^{-50} < 1 -  10^{-30} .\]
\end{proof}


\subsection{Proof of Theorem \ref{thm:effective}}
Let $f(z)$ be a meromorphic function in the complex plane which satisfies
\[ |df|(z) \leq 1 \quad \text{for all $z \in \affine$}. \]
Here we \textit{don't} suppose the condition (\ref{pre-estimate}) is satisfied.

The following argument is similar to the argument in Section \ref{subsection:upper bound}.
First we establish a result on small squares.
\begin{proposition}\label{prop:effective estimate on small squares}
Let $K$ be a square of length of side $2r_0$ in the complex plane.
Then we have 
\[ \frac{1}{\mathrm{vol}(K)}\int_K |df|^2 \, dxdy \leq 1 - \varepsilon . \]
\end{proposition}
\begin{proof}
If $|df|(z) < 1- \varepsilon$ for all $z \in K$, it is obvious that 
\[ \frac{1}{\mathrm{vol}(K)}\int_K |df|^2 \, dxdy \leq (1 - \varepsilon )^2 < 1 - \varepsilon . \]
Hence we can suppose that there is a point $z_0 \in K$ such that $|df|(z_0) \geq 1 - \varepsilon$.
Since the length of side of $K$ is $2r_0$, 
there is a circular sector $D$ of radius $r_0$ and angle $\pi /2$ centered at $z_0$ such that 
$D \subset K$. 
Because we have the Euclidean symmetry on the complex plane and the complex projective
line is also a homogeneous space, we can suppose $z_0 = 0$ and $f(z_0) = 0$.
 Then $f(z)$ satisfies the condition (\ref{pre-estimate}). 
 Applying Proposition \ref{prop:estimate on sector} to this situation, we get
\begin{equation}
 \frac{1}{\mathrm{vol}(D)}\int_D |df|^2 \, dxdy \leq 1 - 10^{-30} . \label{integral over sector}
\end{equation}
From $\mathrm{vol}(K) = 4r_0^2$ and $\mathrm{vol}(D) =\pi r_0^2 /4$,   
\[ \frac{\mathrm{vol}(D)}{\mathrm{vol}(K)} = \frac{\pi}{16} . \]
Using (\ref{integral over sector}) and $|df| \leq 1$, we get
\begin{equation*}
 \begin{split}
 \frac{1}{\mathrm{vol}(K)} \int_K |df|^2 \, dxdy &= 
\frac{\mathrm{vol}(D)}{\mathrm{vol}(K)}\cdot \frac{1}{\mathrm{vol}(D)}\int_D |df|^2 \, dxdy
+ \frac{1}{\mathrm{vol}(K)} \int_{K \setminus D} |df|^2 \, dxdy ,\\
& \leq \frac{\mathrm{vol}(D)}{\mathrm{vol}(K)} (1 - 10^{-30}) 
+ \frac{1}{\mathrm{vol}(K)} (\mathrm{vol}(K) - \mathrm{vol}(D) ) , \\
&= 1 - \frac{\pi}{16} \cdot 10^{-30} , \\
&< 1 - 10^{-100} = 1 - \varepsilon .
  \end{split}
\end{equation*}
\end{proof}
\begin{proof}[\indent\sc Proof of Theorem \ref{thm:effective}]
We prove the theorem by packing squares of length of side $2 r_0$ in the disk 
$\bar{\Delta}(R) = \{ z \in \affine \, |\, |z| \leq R \}$.
If we consider a tiling of the complex plane by squares of length of side $2 r_0$,
 the disk $\bar{\Delta}(R - 2\sqrt{2} r_0)$ is covered by the squares contained in $\bar{\Delta}(R)$.
And we can use the estimate of Proposition \ref{prop:effective estimate on small squares} on each squares
of the tiling.
 
Then the rest of the arguments are the same as the proof of Theorem \ref{thm:upper-bound}.
We omit the details.
\end{proof}


\section{Holomorphic capacity and the Nevanlinna theory}
In this section we prove Theorem \ref{thm:Nevanlinna 2} by using 
the Nevanlinna theory. In Section 6.1 we review some basic facts on the Nevanlinna theory in one variable.
In Section 6.2 we establish preliminary estimates. In Section 6.3 we give the description of 
the holomorphic curves in the complement of hyperplanes and prove Theorem \ref{thm:Nevanlinna 2}.

We use the Fubini-Study metric restricted to $\affine P^n \setminus (H_0 \cup \cdots \cup H_n)$ as the metric on it. 

\subsection{Review of the Nevanlinna theory} \label{subsection:review of the Nevanlinna theory}
We set 
\[ \log^+x = \max ( \log x , 0) \quad \text{for all $x \geq 0$} .\]
For a meromorphic function $f(z)$ in the complex plane, we define the Shimizu-Ahlfors characteristic 
function $T(r, f)$ by 
\begin{equation*}
T(r,f) := \int_1^r \frac{dt}{t} \int_{|z| \leq t} f^* \omega_{FS} \quad \text{for all $r\geq 1$}.
\end{equation*} 
Here $\omega_{FS}$ is the Fubini-Study metric form on the complex projective line defined by 
(\ref{def:Fubini-Study}).
We also define $m(r, f)$ by 
\begin{equation*}
m(r,f) := \frac{1}{2\pi} \int_{|z|= r} \log ^+ |f(z)| d\theta 
= \frac{1}{2\pi} \int_0^{2\pi} \log ^+ |f(r e^{\sqrt{-1}\theta})| d\theta \quad \text{for all $r\geq 1$}.
\end{equation*}
Here $(r, \theta)$ is the polar coordinate in the complex plane.

The following facts are standard in the Nevanlinna theory in one variable.

\begin{fact}\label{fact:Shimizu-Ahlfors}
Let $f(z)$ be a \textit{holomorphic} function in the complex plane.
Then 
\[ T(r,f) = m(r,f) + O(1) , \quad (r \to \infty) . \]
\end{fact}
\begin{fact} \label{fact:characteristic of polynomials}
Let $f(z)$ be a \textit{holomorphic} function in the complex plane. 
$f(z)$ becomes a polynomial if the following condition is satisfied:
\[ \liminf_{r \to \infty} \frac{m(r,f)}{\log r}  < \infty . \]
\end{fact}
\begin{fact}[Nevanlinna's lemma on the logarithmic derivative] \label{fact:logarithmic derivative}
For a meromorphic function $f(z)$ in the complex plane, there are a positive constant $C$ and 
a Lebesgue measurable set $E \subset [1, \infty)$ with a finite measure such that 
\[ m(r, f'/f) \leq C (\log ^+ T(r,f) + \log r ) \quad \text{for all $r \in [1, \infty) \setminus E$} . \]
\end{fact}
\begin{remark}
Actually, the above results are special cases of more general and stronger theorems.
\end{remark} 

We also need the following easy lemma.
\begin{lemma}\label{lem:estimate of exp(polynomial)}
Let $g(z)$ be a polynomial of degree $n \geq 1$. Set $f(z) := e^{g(z)}$. Then 
we have positive constants $r_0$ and  $C$ such that 
\[ m(r, f) \geq C r^n \quad \text{for all $r\geq r_0$} . \]
\end{lemma}
\begin{proof}
Let $g(z) = a_0 z^n + a_1 z^{n-1} + \cdots + a_n$ with $a_0 \neq 0$. Using the rotation of the coordinate, 
we can suppose that $a_0$ is a positive real number. 
 We have a positive constant $C_1$ such that
 \[ |a_1 z^{n-1} + a_2 z^{n-2} + \cdots + a_n | \leq C_1 r^{n-1} , \quad (r \geq 1) . \]
Hence
\begin{equation*}
 \begin{split}
 \mathrm{Re}(g(z)) &= a_0 r^n \cos n\theta + \mathrm{Re}(a_1 z^{n-1} + \cdots + a_n) , \\
 &\geq a_0 r^n \cos n\theta - C_1 r^{n-1} , \quad (r \geq 1).
 \end{split}
\end{equation*}
From $|f(z)| = e^{\mathrm{Re}(g(z))}$,
\[ \log^+ |f(z)| \geq \mathrm{Re}(g(z)) 
\geq a_0 r^n \cos n\theta - C_1 r^{n-1}, \quad (r \geq 1). \]
Then 
\begin{equation*}
 \begin{split}
 \frac{1}{2\pi} \int_{|z|= r} \log ^+ |f(z)| d\theta &\geq \frac{1}{2\pi}\int_0^{\pi /2n} \log^+ 
|f(r e^{\sqrt{-1}\theta})| d\theta ,\\
 &\geq \frac{a_0 r^n}{2\pi} \int_0^{\pi /2n} \cos n\theta d\theta -  \frac{C_1}{4n}\, r^{n-1} , \\
 &= \frac{a_0 r^n}{2\pi n} - \frac{C_1}{4n}\, r^{n-1} , \quad (r \geq 1).
 \end{split}
\end{equation*}
Thus we have positive constants $r_0$ and $C$ such that 
\[ m(r, f) \geq C r^n , \quad (r \geq r_0). \]
\end{proof}


\subsection{Preliminary estimates}

For a meromorphic function $f(z)$ in the complex plane, 
\[ f^* \omega_{FS} = |df|^2 \, dxdy = \frac{1}{\pi} \frac{|f'(z)|^2}{(1 + |f(z)|^2)^2} dxdy .\]

\begin{lemma}\label{lem:estimate of exp(az+b)}
Let $a$ and $b$ be complex numbers and set $f(z) := e^{az + b}$. 
Then we have a positive constant $C$ such that 
\[ \int_{|z| \leq r} |df|^2\, dxdy \leq C r \quad \text{for all $r \geq 0$}. \]
\end{lemma}
\begin{proof}
If $a= 0$, the statement is trivial. Hence we can suppose that $a \neq 0$. 
Using the rotation of the coordinate, we can also suppose that $a$ is a positive real number. 
Set $\beta := \mathrm{Re}(b)$. Then we have
\[ |df|^2(z) = \frac{a^2}{\pi}\frac{e^{2ax+ 2\beta}}{(1+ e^{2ax+ 2\beta})^2} .\]
Thus
\begin{equation*}
 \begin{split}
 \int_{|z| \leq r} |df|^2\, dxdy &\leq \int_{-r}^r dy \int_{-\infty}^{\infty} dx \,
 \frac{a^2}{\pi}\frac{e^{2ax+ 2\beta}}{(1+ e^{2ax+ 2\beta})^2} , \\
 &= \frac{2 a^2 r}{\pi}\int_{-\infty}^{\infty} \frac{e^{2ax+ 2\beta}}{(1+ e^{2ax+ 2\beta})^2} \, dx ,\\ 
 &= \frac{ar}{\pi}.
 \end{split}
\end{equation*}
\end{proof}

For a Lebesgue measurable set $E$ in $\mathbb{R}$, we denote its Lebesgue measure by $|E|$.
\begin{lemma}\label{lem:estimate of exp(az^2+bz+c)}
Let $a,b,c$ be complex numbers with $a\neq 0$, and set $f(z) := e^{az^2+bz+c}$.
Then, for any positive number $\varepsilon$, we have a open subset $E$ in $[0, 2\pi]$ 
with $|E| < \varepsilon$ such that 
\[ \int_{[0, 2\pi]\setminus E} d\theta \int_0^{\infty} |df|^2 (r e^{\sqrt{-1} \theta}) \, r dr  < \infty .\]
\end{lemma}
\begin{proof}
Using the rotation of the coordinate, we can suppose that $a$ is a positive real number.
From $f'(z) = (2az + b)f(z)$, we have
\begin{equation*}
 \begin{split}
 |df|(z) &= \frac{1}{\sqrt{\pi}}\frac{|2az + b||f(z)|}{1 + |f(z)|^2}, \\
 |f(z)| &= \exp (\mathrm{Re}(az^2 + bz +c)).
 \end{split}
\end{equation*}

We define $E \subset [0, 2\pi]$ by 
\begin{equation*}
 \begin{split}
 E := &\left( \frac{\pi}{4}-\frac{\delta}{2}, \frac{\pi}{4}+\frac{\delta}{2}\right) \cup
        \left( \frac{3 \pi}{4}-\frac{\delta}{2}, \frac{3 \pi}{4}+ \frac{\delta}{2} \right)  \\
      &\cup  \left( \frac{5 \pi}{4}-\frac{\delta}{2}, \frac{5 \pi}{4} + \frac{\delta}{2}\right) \cup
       \left( \frac{7 \pi}{4}-\frac{\delta}{2}, \frac{7 \pi}{4}+ \frac{\delta}{2}\right) .
 \end{split}
\end{equation*}
Here $\delta$ is a sufficiently small positive number such that $|E| = 4\delta < \varepsilon$.
Then, we have 
\[ |\cos 2\theta | \geq \sin \delta \quad \text{for all $\theta \in [0, 2\pi]\setminus E$}. \]
\noindent
(i)
If $\cos 2\theta \geq \sin \delta$, we have
\begin{equation*}
 \begin{split}
 \mathrm{Re}(az^2 + bz + c) &= a r^2 \cos 2\theta + \mathrm{Re}(bz + c) , \\
 &\geq a (\sin \delta )r^2 - |b|r - |c|.
 \end{split}
\end{equation*}
Hence
\begin{equation*}
 \begin{split}
 |df|(z) &\leq \frac{1}{\sqrt{\pi}}\frac{|2az + b||f(z)|}{|f(z)|^2} ,\\
 &\leq \frac{1}{\sqrt{\pi}} (2ar + |b|) \exp (- \mathrm{Re}(az^2 + bz + c)) ,\\
 &\leq \frac{1}{\sqrt{\pi}} (2ar + |b|) \exp ( - a (\sin \delta )r^2 + |b|r + |c|).
 \end{split}
\end{equation*} 
\noindent
(ii)
If $\cos 2\theta \leq - \sin \delta$, we have 
\begin{equation*}
 \begin{split}
 \mathrm{Re}(az^2 + bz + c) &= ar^2 \cos 2\theta + \mathrm{Re}(bz + c) ,\\
 &\leq -a (\sin \delta )r^2 + |b|r + |c|.
 \end{split}
\end{equation*}
Hence
\begin{equation*}
 \begin{split}
 |df|(z) &\leq \frac{1}{\sqrt{\pi}}|2az + b| |f(z)| ,\\
 &\leq \frac{1}{\sqrt{\pi}}  (2ar + |b|) \exp (\mathrm{Re}(az^2 + bz + c)) ,\\
 &\leq \frac{1}{\sqrt{\pi}} (2ar + |b|) \exp (-a (\sin \delta )r^2 + |b|r + |c|).
 \end{split}
\end{equation*} 

Therefore we get 
\[ |df|(r e^{\sqrt{-1}\theta}) \leq 
\frac{1}{\sqrt{\pi}} (2ar + |b|) \exp (-a (\sin \delta )r^2 + |b|r + |c|) \quad 
\text{if $\theta \in [0, 2\pi]\setminus E$}. \]
Thus we have a constant $C$ independent of $\theta$ such that 
\[ \int_0^{\infty} |df|^2 (r e^{\sqrt{-1}\theta})\, rdr \leq C 
\quad \text{for all $\theta \in [0, 2\pi]\setminus E$}. \]
It follows that 
\[ \int_{[0, 2\pi]\setminus E} d\theta \int_0^{\infty} |df|^2 (r e^{\sqrt{-1} \theta}) \, r dr 
 \leq 2\pi C < \infty . \] 
\end{proof}


\subsection{Holomorphic curves in the complement of hyperplanes}

We define the $n+1$ hyperplanes $P_0, P_1, \cdots , P_n$ in $\affine P^n$ by 
\[ P_i :  z_i = 0, \quad (0 \leq i \leq n) . \]
Here $[z_0, z_1, \cdots , z_n]$ is the homogeneous coordinate of $\affine P^n$.
Then
\begin{equation}\label{the complement of hyperplanes}
 \affine P^n \setminus (P_0 \cup \cdots \cup P_n) = \{\, [1:z_1:z_2:\cdots :z_n]\,|\, z_i \neq 0 , 
\, (1 \leq i \leq n) \} \cong (\affine \setminus \{0 \})^n .
\end{equation}

The following proposition is proved in [BD, Appendice].
\begin{proposition}
Let $f: \affine \to \affine P^n \setminus (P_0\cup \cdots \cup P_n)$ be a holomorphic map with 
$\norm{df} < \infty$. Then, there are complex numbers $a_i$ and $b_i$, $(1\leq i \leq n)$, such that
\[ f(z) = [\, 1: \exp (a_1 z + b_1): \exp (a_2 z + b_2): \cdots : \exp (a_n z + b_n)\,] . \]
\end{proposition}
\begin{proof}
\footnote{This proof is less conceptual than the proof in [BD]. 
But the author think that our proof also contains some interesting arguments.}
From (\ref{the complement of hyperplanes}), we have holomorphic maps 
$f_i:\affine \to \affine \setminus \{ 0 \}$, $(1\leq i \leq n)$, such that 
\[ f(z) = [\, 1: f_1(z): f_2(z): \cdots : f_n(z)\, ]. \]
Since $\exp : \affine \to \affine \setminus \{ 0 \}$ is the universal covering, there are holomorphic 
functions $g_i(z)$ in the complex plane such that $f_i(z) = \exp (g_i(z))$, $(1 \leq i \leq n)$. 
We will prove that all $g_i(z)$ are linear functions.
The proof falls into three steps. First we prove that $g_i(z)$ are polynomials.
Next we show $\deg (g_i(z)) \leq 2$.
 In the last step we prove $\deg (g_i(z)) \leq 1$.
The arguments in the first and second steps are standard. The last step is a little tricky.

Set $m := \norm{df}$. Here $|df|^2 = \frac{1}{4\pi} \Delta \log (1 + \sum |f_i|^2)$. 
Using Jensen's formula, we get 
\begin{equation}\label{estimate:Jensen's formula}
 \begin{split}
 \int_1^r \frac{dt}{t} \int_{|z| \leq t} \frac{1}{4\pi} &\Delta \log \left(1 + \sum_i |f_i|^2\right) \, dxdy \\
 &= \frac{1}{4\pi} \int_{|z| = r} \log \left( 1 + \sum_i |f_i|^2\right) d\theta 
    -  \frac{1}{4\pi} \int_{|z| = 1} \log \left( 1 + \sum_i |f_i|^2\right) d\theta . 
 \end{split}
\end{equation}
Hence 
\begin{equation*}
 \begin{split}
 \frac{1}{4\pi} \int_{|z| = r} \log \left( 1 + \sum |f_i|^2\right) d\theta
 &=  \int_1^r \frac{dt}{t} \int_{|z| \leq t} |df|^2 \, dxdy  + \mathrm{const}, \\
 &\leq \int_1^r \frac{dt}{t} \int_{|z| \leq t} m^2 \, dxdy  + \mathrm{const},  \\
 &\leq \frac{1}{2}m^2 \pi r^2 + \mathrm{const} .
 \end{split}
\end{equation*}
Since we have
\[ \log ^+ |f_i| = \frac{1}{2} \log ^+|f_i|^2 \leq \frac{1}{2} \log (1 + |f_1|^2 + \cdots + |f_n|^2) ,\]
it follows that 
\begin{equation}\label{estimate of m(r,f_i)}
 \begin{split}
 m(r, f_i) &\leq \frac{1}{4\pi} \int_{|z| = r} \log \left( 1 + \sum_j |f_j|^2\right) d\theta , \\
 & \leq \frac{1}{2}m^2 \pi r^2 + \mathrm{const} .
 \end{split}
\end{equation}
From Fact \ref{fact:Shimizu-Ahlfors}, we get
\[ T(r, f_i) \leq \frac{1}{2}m^2 \pi r^2 + \mathrm{const}' . \]
Fact \ref{fact:logarithmic derivative} gives 
\begin{equation*}
 \begin{split}
 m(r, g_i') = m(r, f_i'/f_i) &\leq C (\log ^+ T(r, f_i) + \log r) ,\\
&\leq  \mathrm{const}\cdot \log r + \mathrm{const} \quad \text{for all $r \in [1, \infty) \setminus E$}. 
 \end{split}
\end{equation*}
Here $E$ is a Lebesgue measurable set in $[1, \infty)$ with a finite measure.
It follows that 
\[ \liminf_{r \to \infty} \frac{m(r, g_i')}{\log r}  < \infty . \]
From Fact \ref{fact:characteristic of polynomials}, this shows that $g_i'(z)$ is a polynomial.
Hence $g_i(z)$ is also a polynomial.
 
Next we will prove $\deg (g_i(z)) \leq 2$. Suppose, for example, $\deg (g_1(z)) \geq 3$.
From Lemma \ref{lem:estimate of exp(polynomial)} 
and the estimate of $m(r, f_1)$ in (\ref{estimate of m(r,f_i)}),
 we have a positive constant $C_1$ such that
\[ C_1 r^3 \leq m(r,f_1) \leq \frac{1}{2}m^2 \pi r^2 + \mathrm{const}, \quad (r \gg 0) . \]
This is obviously impossible. Hence $\deg (g_i(z)) \leq 2$, $(1 \leq i \leq n)$.

Finally we will prove $\deg (g_i(z)) \leq 1$. Suppose, for example, $\deg (g_1(z)) = 2$.
From Lemma \ref{lem:estimate of exp(polynomial)}, we have positive constants $r_0$ and $C_2$ such that
\[ m(r, f_1) \geq C_2 r^2, \quad (r \geq r_0). \]
From (\ref{estimate:Jensen's formula}) and (\ref{estimate of m(r,f_i)}), 
\begin{equation}\label{estimate of int |df|^2}
 \begin{split}
C_2 r^2 \leq m(r, f_1) &\leq \frac{1}{4\pi} \int_{|z| = r} \log \left( 1 + \sum_i |f_i|^2\right) d\theta ,\\
&= \int_1^r \frac{dt}{t} \int_{|z| \leq t} |df|^2 \, dxdy  + \mathrm{const}, 
\quad (r\geq r_0).
 \end{split}
\end{equation}
From the definition of the Fubini-Study metric (\ref{def:Fubini-Study}), 
\begin{equation} \label{estimate of |df|^2 by f_i and f_i/f_j}
 \begin{split}
 |df|^2 &= \frac{1}{\pi}\left[ \frac{\sum_i |f'_i|^2}{(1 + \sum_i |f_i|^2)^2} 
+ \frac{\sum_{i<j} |g'_i - g'_j|^2 |f_i|^2 |f_j|^2}{(1 + \sum_i |f_i|^2)^2} \right] , \\  
 &\leq \frac{1}{\pi} \left[ \sum_i \frac{|f'_i|^2}{(1 + |f_i|^2)^2} 
 + \sum_{i<j} \frac{|g'_i - g'_j|^2 |f_i|^2 |f_j|^2}{(|f_i|^2 + |f_j|^2)^2}\right] , \\
 &= \frac{1}{\pi} \left[ \sum_i \frac{|f'_i|^2}{(1 + |f_i|^2)^2} 
 + \sum_{i<j} \frac{|(f_i/f_j)'|^2}{(1 + |f_i/f_j|^2)^2} \right] .
 \end{split}
\end{equation}
Because $f_i(z) = \exp (g_i(z))$, $f_i(z)/f_j(z) = \exp (g_i(z) - g_j(z))$, and 
the degrees of $g_i(z)$ and $g_i(z) - g_j(z)$ are at most two,
we can apply Lemma \ref{lem:estimate of exp(az+b)} and Lemma \ref{lem:estimate of exp(az^2+bz+c)}
to holomorphic functions $f_i(z)$ and $f_i(z)/ f_j(z)$.
Hence there are positive constants $C_3$, $C_4$ and
 a open set $F \subset [0, 2\pi]$ with $|F| < 2C_2/m^2$ such that 
\[ \int_{[0, 2\pi]\setminus F} d\theta \int_0^t |df|^2(r e^{\sqrt{-1}\theta}) \, rdr
\leq C_3 t + C_4 . \]
Then we have 
\begin{equation*}
 \begin{split}
 \int_{|z| \leq t} |df|^2 \, dxdy  &= 
 \int_{F} d\theta \int_0^t |df|^2(r e^{\sqrt{-1}\theta}) \, rdr
 + \int_{[0, 2\pi]\setminus F} d\theta \int_0^t |df|^2(r e^{\sqrt{-1}\theta}) \, rdr ,\\
 &\leq \int_{F} d\theta \int_0^t m^2 \, rdr + C_3 t + C_4 ,\\
 &= \frac{m^2 |F|}{2}t^2 + C_3 t + C_4 , \\
 &\leq C_2 t^2 + C_3 t + C_4 .
 \end{split}
\end{equation*}
Substituting this into (\ref{estimate of int |df|^2}), we get
\[ C_2\, r^2 \leq \frac{C_2}{2} \, r^2 + C_3 r + C_4 \log r + \mathrm{const} , \quad (r\geq r_0). \]
This is impossible. Thus we conclude that all $g_i(z)$ are linear functions.
\end{proof}

\begin{proposition}\label{prop:holomorphic capaity of the complement of hyperplanes}
Let $f: \affine \to \affine P^n \setminus (P_0 \cup \cdots \cup P_n)$ be a holomorphic map with 
$\norm{df} < \infty$. Then 
\[ \limsup_{R \to \infty} \frac{1}{\pi R^2} \int_{|z| \leq R} |df|^2 \, dxdy = 0 . \]
\end{proposition}
\begin{proof}
From (\ref{estimate of |df|^2 by f_i and f_i/f_j}),
\[ |df|^2 \leq \frac{1}{\pi} \left[ \sum_i \frac{|f'_i|^2}{(1 + |f_i|^2)^2} 
 + \sum_{i<j} \frac{|(f_i/f_j)'|^2}{(1 + |f_i/f_j|^2)^2} \right] . \]
Since $f_i(z) = \exp (g_i(z))$ and $f_i(z)/f_j(z) = \exp (g_i(z) - g_j(z))$ with linear functions
$g_i(z)$ and $g_i(z) - g_j(z)$,
we can apply Lemma \ref{lem:estimate of exp(az+b)} and get 
\[ \int_{|z|\leq R} |df|^2 \, dxdy \leq C R. \]
Here $C$ is a positive constant. Thus
\[ \limsup_{R \to \infty} \frac{1}{\pi R^2} \int_{|z| \leq R} |df|^2 \, dxdy \leq 
\lim_{R \to \infty} \frac{C}{\pi R} = 0. \]
\end{proof}

\begin{proof}[\indent\sc Proof of Theorem \ref{thm:Nevanlinna 2}]
Using the defining equations of $H_i$ in (\ref{defining equations of hyperplanes}),
 we define a biholomorphic map 
$A: \affine P^n \to \affine P^n$ by 
\[ A([z_0:z_1:\cdots :z_n]) := 
\left[ \sum_{j} a_{0j}z_j : \sum_{j} a_{1j}z_j: \cdots :\sum_{j} a_{nj}z_j \right] . \]
$A$ gives a biholomorphic map from $\affine P^n \setminus (H_0\cup \cdots \cup H_n)$ to 
$\affine P^n \setminus (P_0\cup \cdots \cup P_n)$.
Since $A$ and $A^{-1}$ are defined on the \textit{compact} set $\affine P^n$, 
we have a constant $C$ such that 
\[ |d A (u)| \leq C |u| \quad \text{and} \quad
|d A^{-1} (u)| \leq C |u| \quad \text{for all $u \in T\affine P^n$}. \]

Let $f: \affine \to \affine P^n \setminus (H_0\cup \cdots \cup H_n)$ be a holomorphic map
 with $\norm{df} \leq 1$. Then $Af$ is a holomorphic map from $\affine$ to
 $\affine P^n \setminus (P_0\cup \cdots \cup P_n)$ with $\norm{d(Af)} \leq C <\infty$.
Since $f = A^{-1}\circ Af$, Proposition \ref{prop:holomorphic capaity of the complement of hyperplanes}
gives
\[ \rho (f) = \limsup_{R \to \infty} \frac{1}{\pi R^2} \int_{|z|\leq R} |df|^2 \, dxdy
  \leq \limsup_{R \to \infty} \frac{C^2}{\pi R^2} \int_{|z|\leq R} |d(Af)|^2\, dxdy = 0.   \] 
Thus we conclude that 
\[ \rho (\affine P^n \setminus (H_0 \cup \cdots \cup H_n)) = 0 .\]
\end{proof}

\vspace{10mm}

Masaki Tsukamoto, Department of Mathematics, Faculty of Science, Kyoto University, Kyoto 606-8502, Japan

\textit{E-mail address}: \texttt{tukamoto@math.kyoto-u.ac.jp}

\end{document}